\documentclass[11pt]{article}
\usepackage{amsmath,amssymb}
\usepackage[mathscr]{euscript}

\usepackage[parfill]{parskip} %Here, we delete indent and space between paragraphs
\newcommand{\proof}{\vspace{0pt}{\bfseries Proof.\ }}
\newcommand{\QED}{{\unskip\nobreak\hfil\penalty50\hskip1em\hbox{}\nobreak
   \hfil \ensuremath{\Box}\parfillskip=0pt \par}}
\newtheorem{Theorem}{Theorem}[section]
\newtheorem{Lemma}[Theorem]{Lemma}
\newtheorem{Corollary}[Theorem]{Corollary}
\newtheorem{Proposition}[Theorem]{Proposition}
\newtheorem{Remark}[Theorem]{Remark}

\newtheorem{Definition}[Theorem]{Definition}

\begin{document}

\begin{center}
{\large\bf Complete Intersection Toric Ideals of Oriented Graphs and Chorded-Theta Subgraphs
}

\medskip

I. Gitler, E. Reyes, J.A. Vega \\ 
{\small Departamento de Matem\'aticas}\vspace{-1mm}\\ 
{\small Centro de Investigaci\'on y de Estudios Avanzados del
IPN}\vspace{-1mm}\\   
{\small Apartado Postal 14--740}\vspace{-1mm}\\ 
{\small 07000 M\'exico City, D.F.}\vspace{-1mm}\\ 
{\small e-mail: {\tt ereyes@math.cinvestav.mx}}\vspace{4mm}
\end{center}
\date{}

\begin{abstract} 
\noindent
Let $G=(V,E)$ be a finite, simple graph. We consider for each oriented graph $G_{\cal O}$ associated to an orientation ${\cal O}$ of the edges of $G$, the toric ideal $P_{G_{\cal O}}$. In this paper we study those graphs with the property that $P_{G_{\cal O}}$ is a binomial complete intersection, for all ${\cal O}$. These graphs are called $\text{CI}{\cal O}$ graphs. We prove that these graphs can be constructed recursively as clique-sums of cycles and/or complete graphs. We introduce the chorded-theta subgraphs and their transversal triangles. Also we establish that the $\text{CI}{\cal O}$ graphs are determined by the property that each chorded-theta has a transversal triangle. As a consequence, we obtain that the tournaments hold this property. Finally we explicitly give the minimal forbidden induced subgraphs that characterize these graphs, these families of graphs are: prisms, pyramids, thetas and a particular family of wheels that we call $\theta-$partial wheels.

\medskip

\end{abstract}

\section{Introduction}\label{Int}

Let $G=(V,E)$ be a graph with $ \vert E({G}) \vert =q$ and whose vertex set is given by $V(G)=\{x_1,\dotsc,x_n\}$. An orientation ${\cal O}$ of the edges of $G$ is an assignment of a direction to each edge of $G$. Let $D=G_{\cal O}$ denote the oriented graph associated to an orientation ${\cal O}$ of the edges of $G$. In particular if $G$ is a complete graph, $G_{\cal O}$ is called a {\it tournament}. To each oriented edge $e=(x_i,x_j)$ of $D$, we associate the vector $v_e \in \{0,1,-1\}^n$ defined as follows: the $i${\it th} entry 
is $-1$, the $j${\it th} entry is $1$, and the remaining entries are zero. The {\it incidence matrix\/} $A_{D}$ of $D$ is the $n\times q$ matrix whose columns are the vectors of the form $v_e$, with $e$ an edge of $D$. The set of column vectors of $A_D$ will be denoted by ${\cal A}=\{v_1,\dotsc,v_q\}$.  
\\
\\
Consider the {\it edge subring\/} $k[D]:=k[x^{v_1},\dotsc,x^{v_q}]\subset k[x_1^{\pm 1},\dotsc,x_n^{\pm 1}]$ of $D$, where $v_i=(v_i^1,\dotsc,v_i^n)$ and $x^{v_i}=x_1^{v_i^1}\dotsm x_n^{v_i^n}$. Let $E(D)=\{t_1,\dotsc,t_q\}$ be the edge set of $D$. There is an epimorphism of $k-$algebras given by
\[
\varphi\colon k[t_1,\dotsc,t_q] \longrightarrow k[D], \ \ \ 
\mbox{ where } t_i\longmapsto x^{v_i},
\]
and $k[t_1,\dotsc,t_q]$ is a polynomial ring. The kernel of $\varphi$, denoted by $P_D$, is called the {\it toric ideal\/} of $D$. This ideal has been thoroughly studied in \cite{GRV,Ring,morf}. The toric ideal $P_D$ is a prime ideal of height $q-n+r$, where $r$ is the number of connected components of $D$. $P_D$ is generated by binomials and $k[D]$ is a normal domain (see \cite{Ring,monalg}). Thus any minimal generating set of $P_D$ must have at least $q-n+r$ elements, by the Principal Ideal Theorem (see \cite{AM}). If $P_D$ can be generated exactly by $q-n+r$ polynomials it is called a {\it complete intersection}. If these polynomials are binomials then $P_D$ is called a {\it binomial complete intersection}. In Section \ref{Toric-Orient} we study the binomials of $P_D$. In particular, we prove that the primitive binomials of $P_D$ are the binomials associated to cycles. We also recover the result given in \cite{Ring} that $P_D$ is generated by binomials associated to cycles without chords. If two graphs $G$ and $H$ each contain cliques of equal size, the {\it clique-sum\/} of $G$ and $H$ is formed from their disjoint union by identifying pairs of vertices in these two cliques to form a single shared clique.
\\
\\
It is known  \cite{Ring,morf} that any graph has at least one acyclic orientation such that the corresponding 
toric ideal is a binomial complete intersection. It is natural to ask what is the class of graphs with the property that $P_{G_{\cal O}}$ is generated by exactly $q-n+r$ binomials for every orientation ${\cal O}$ of $G$. For ease of exposition these graphs are called $\text{CI}{\cal O}$ graphs. In the bipartite case a graph $G$ is in this class if and only if  $G$ is a ring graph (see \cite{Ring}). Ring graphs are introduced in \cite{GRV}, they are those graphs for which their non trivial blocks (blocks that are not vertices or bridges) can be constructed as $2-$clique-sums of cycles. In Section \ref{CIO-graphs} we show that the $\text{CI}{\cal O}$ property is closed under induced subgraphs and we prove that $G$ is a $\text{CI}{\cal O}$ graph if and only if its blocks are $\text{CI}{\cal O}$ graphs. 
\\
\\
In Section \ref{Ch-thetas} we introduce the family of chorded-theta subgraphs and their transversal triangles. In this section we study the graphs $G$ that satisfy: for every chorded-theta of $G$ there exists a transversal triangle ($\forall \theta \exists \Delta-$property). We call this type of graphs theta-ring graphs.
In Theorem \ref{For-all-Theta-exists-delta} we prove that these graphs can be constructed as clique-sums of cycles and complete graphs. Furthermore we show that they can be obtained as $0,1,2-$clique-sums  of chordal graphs and cycles. Theta-ring graphs are closed under induced subgraphs and we obtain the corresponding set of minimal forbidden induced subgraphs that characterize this family. These minimal forbidden induced subgraphs are: prisms, pyramids, thetas and $\theta-$partial wheels (partial wheels that are chorded-theta graphs).  
\\
\\ 
Our main result in Section \ref{CIO-Theta} is Theorem \ref{For-all-Theta-exists-deltaCIO}: $P_{G_{\cal O}}$ is a binomial complete intersection for every orientation ${\cal O}$ of $G$ if and only if $G$ is a theta-ring graph, i.e., the $\text{CI}{\cal O}$ property is equivalent to the $\forall \theta \exists \Delta-$property. As a consequence of this result, we obtain that $P_{G_{\cal O}}$ is a binomial complete intersection if $G_{\cal O}$ is a tournament. Using the results given in Section \ref{Ch-thetas} and Section \ref{CIO-Theta} we obtain the following equivalences: 
\[
\begin{array}{cccccc}
  \begin{array}{c} G \mbox{ is a }\\ 
        \mbox{ $\text{CI}{\cal O}$ graph }
  \end{array}
&\Leftrightarrow & 
\begin{array}{c} G \mbox{ is a theta-ring graph}\\
                       (G \mbox{ has the }\forall \theta \exists \Delta \mbox{-property})
                                 \end{array}& 
%                                 \Leftrightarrow 
%& \begin{array}{c}\mbox{For every chorded-theta} \\ 
%                        \mbox{of $G$ there exists a}\\
%                        \mbox{transversal triangle}\\
%                         \end{array} &
\Leftrightarrow &\begin{array}{c}\mbox{$G$ is obtained by clique-sums}\\ 
                        %\mbox{  }\\  
                        \mbox{of complete graphs and cycles}                     
                   \end{array}&\\     
\end{array}
\] 
\\
\\
In particular, we recover that ring graphs are $\text{CI}{\cal O}$ graphs and that the converse also holds in the bipartite case. 
\\
\\
The paper is essentially self contained. For unexplained terminology and notation on graph theory and toric ideals we refer to \cite{diestel,Ring,gitler-vila,monalg}. Some references for toric ideals associated to graphs (without orientation) are \cite{GRV,Ring,Katsabekis,katzman,Reyes-Tatakis-Thoma,aron-jac,Tatakis-Thoma,Vi3}.

\section{Toric ideals of oriented graphs}\label{Toric-Orient}

Let $G$ be a graph and consider the oriented graph $D=G_{\cal O}$ associated to an orientation ${\cal O}$ of the edges of $G$. Let $\alpha=(\alpha_1,\dotsc,\alpha_q)\in\mathbb{R}^q$, the support of $\alpha$ is the set $\text{supp}(\alpha)=\{i \mid \alpha_i \neq 0\}$. Furthermore $\alpha=\alpha_+-\alpha_-$, where $\alpha_+$ and $\alpha_-$ are two non negative vectors with disjoint support. We denote by $P_D$ the kernel of the epimorphism of $k-$algebras
\[ 
\varphi\colon k[t_1,\dotsc,t_q] \longrightarrow k[D], \ \ \ 
\mbox{ where } t_i\longmapsto x^{v_i}.
 \]
On the other hand we have the following linear map 
\[ 
\psi\colon {\mathbb Z}^q \longrightarrow {\mathbb Z}^n, \ \ \ 
\mbox{ where } e_i\longmapsto v_i,
\]
i.e., $\psi$ is the linear map associated to the incidence matrix $A_D$ of $D$. If $0\neq\alpha\in {\mathbb Z}^q$ we associate to it the binomial $t_\alpha=t^{\alpha_+}-t^{\alpha_-}$. Notice that $\alpha \in$ $\ker(\psi)$ if and only if $t_\alpha\in P_D$. Given a cycle $C=(x_{j_1},x_{j_2},\dotsc,x_{j_k},x_{j_{k+1}}=x_{j_1})$ of $D$, we split the edge set of $C$ in two disjoint sets $C_+$ and $C_-$, where $C_+$ is the set of clockwise oriented edges, i.e., $C_+=\{t_s \in E(C) \mid \mbox{There exists $i$ such that } t_s=(x_{j_i},x_{j_{i+1}}) \}$ and $C_-=C\setminus C_+$. Associated to $C$ we have the binomial 
\[
t_C=\prod_{t_i\in C_+}t_i-\prod_{t_i\in C_-}t_i.
\]
This binomial belongs to $P_D$. If $C_+=\emptyset$ or $C_-=\emptyset$ we set $\prod\limits_{t_i\in C_+}t_i=1$ or $\prod\limits_{t_i\in C_-}t_i=1$.
\\
\\
In this section we study some properties of the binomials in $P_D$. In Theorem \ref{bincycle} we prove that if $b$ is a binomial in $P_D$ then, there exists a binomial $b'$ associated to a cycle whose monomials divide the monomials of $b$. Using this result we prove that the primitive binomials in $P_D$ are the ones associated with cycles and we recover the result in \cite{Ring} that $P_D$ is generated by the binomials corresponding to cycles without chords.

\begin{Lemma}\label{binedge}
Let $t^{\alpha}-t^{\beta}$ be a binomial in $P_D$ with $\gcd(t^{\alpha},t^{\beta})=1$. If $i\in\text{supp}(\alpha)$ and $t_i=(x_{j_1},x_{j_2})$ then there exist edges $t_{i_1}$ and $t_{i_2}$ both different from $t_i$ such that $x_{j_1} \in V(t_{i_1})$ and $x_{j_2} \in V(t_{i_2})$. Furthermore if $t_{i_2}=(x_{j_2},x_{j_3})$, then $t_{i_2} \mid t^{\alpha}$ and if $t_{i_2}=(x_{j_3},x_{j_2})$ then $t_{i_2} \mid t^{\beta}$.
\end{Lemma}
\proof
Since $t^{\alpha}-t^{\beta} \in P_D$ then 
\[
\varphi(t_1^{\alpha_1}) \dotsm \varphi(t_q^{\alpha_q})=\varphi(t^{\alpha})=\frac{m_1}{m_2}=\varphi(t^{\beta})=\varphi(t_1^{\beta_1}) \dotsm\varphi(t_q^{\beta_q})
\] 
where $\alpha=(\alpha_1,\dotsc,\alpha_q)$ and $\beta=(\beta_1,\dotsc,\beta_q)$. We can assume that $m_1$ and $m_2$ are monomials with $\gcd(m_1,m_2)=1$. Since $t_i=(x_{j_1},x_{j_2})$, we have that $\varphi(t_i^{\alpha_i})=\frac{x_{j_2}^{\alpha_i}}{x_{j_1}^{\alpha_i}}$ and $\alpha_i \neq 0$. If $x_{j_2} \mid m_1$ then there exists an edge $t_{i_2}$ such that ${\beta}_{i_2} \neq 0$ and $\varphi(t_{i_2})= \frac{x_{j_2}}{x_{j_3}}$ for some $x_{j_3}$. Thus $t_{i_2}=(x_{j_3},x_{j_2})$ and $t_{i_2} \mid t^{\beta}$. Furthermore $t_i \neq t_{i_2}$ because $\gcd(t^{\alpha},t^{\beta})=1$. Now, if $x_{j_2}$ does not divide $m_1$ then there exists an edge $t_{i_2}$ such that $\varphi(t_{i_2})=\frac{x_{j_4}}{x_{j_2}}$ with $\alpha_{i_2} \neq 0$ and $i \neq i_2$. Hence, $t_{i_2}=(x_{j_2},x_{j_4})$ and $t_{i_2} \mid t^{\alpha}$, so $x_{j_2} \in V(t_{i_2})$ and 
$t_{i_2} \neq t_i$.
\\
In a similar form we prove that there exists an edge $t_{i_1}$ different of $t_i$ such that $x_{j_1} \in V(t_{i_1})$.
\QED

\begin{Theorem}\label{bincycle}
If $0 \neq t^{\alpha}-t^{\beta} \in P_D$ then there exists a cycle $C$ such that $t_C=t^{\alpha'}-t^{\beta'}$, $t^{\alpha'} \mid t^{\alpha}$ and $t^{\beta'} \mid t^{\beta}$.
\end{Theorem}
\proof 
We can suppose that $\gcd(t^{\alpha},t^{\beta})=1$, otherwise we would have that $\gcd(t^{\alpha},t^{\beta})=m$. Thus, $t^{\alpha}=mt^{\gamma}$, $t^{\beta}=mt^{w}$, $\gcd(t^{\gamma},t^{w})=1$ and $t^{\gamma}-t^{w} \in P_D$. We set $\alpha=(\alpha_1,\dotsc,\alpha_{q})$ and $\beta=(\beta_1,\dotsc,\beta_q)$. Since $t^{\alpha}-t^{\beta} \neq 0$, then there exists $i$ such that $\alpha_i \neq 0$ or $\beta_i \neq 0$. We can suppose that $\alpha_i \neq 0$. By Lemma \ref{binedge} if $x_{j_1} \in V(t_i)$ there exists an edge $t_{i_1}$ different to $t_i$ such that $x_{j_1} \in V(t_i)$ and $t_{i_1} \mid t^{\alpha}t^{\beta}$. Now if $V(t_{i_1})=\{x_{j_1},x_{j_2}\}$ then by Lemma \ref{binedge} there exists an edge $t_{i_2}$ different from $t_{i_1}$ such that $x_{j_2} \in V(t_{i_2})$ and $t_{i_2} \mid t^{\alpha}t^{\beta}$. We continue with this process and we obtain vertices $x_{j_1},x_{j_2}, \dotsc$ such that $\{x_{j_{k}},x_{j_{k+1}}\}\in E(G)$. We take $l= \max \{i \mid x_{j_i} \neq x_{j_k} \mbox{ for } k< i\}$, then $x_{j_{l+1}}=x_{j_s}$ for some $s \leq l$. Therefore $C=(x_{j_s},x_{j_{s+1}},\dotsc,x_{j_l},x_{j_{l+1}}=x_{j_s})$ is a cycle in $G$ and by construction if 
$t_C=t^{\alpha'}-t^{\beta'}$ then $t^{\alpha'}t^{\beta'} \mid t^{\alpha}t^{\beta}$. But, by Lemma \ref{binedge} we have that $t^{\alpha'} \mid t^{\alpha}$ and $t^{\beta'} \mid t^{\beta}$. 
\QED

\begin{Definition}
Let $C=(y_1,y_2, \dotsc,y_k,y_{k+1}=y_1)$ be a cycle of $G$. $C$ is called an oriented cycle in $D=G_{\cal O}$ if $(y_i,y_{i+1}) \in E(D)$ for $i=1,\dotsc,k$ (i.e. $C_+=\emptyset$ or $C_-=\emptyset$). Similarly a path $L=(y_1,\dotsc,y_s)$ of $G$ is called an oriented path in $D$ if $(y_i,y_{i+1}) \in E(D)$ for $i=1,\dotsc,s-1$.
If $L$ is a path, we denote by $L^{o}$ the path $(y_2,\dotsc,y_{s-1})$, if $s=2$ then $L^{o}=\emptyset$.
\end{Definition}

\begin{Corollary}\label{oriencycle}
If $0 \neq 1 - t^{\beta} \in P_D$ then there exists an oriented cycle $C$ in $D$.
\end{Corollary}
\proof 
By Theorem \ref{bincycle}, there is a cycle $C$ such that $t_C=t^{\alpha'}-t^{\beta'}$, $t^{\alpha'} \mid 1$ and $t^{\beta'} \mid t^{\beta}$. Then $t^{\alpha'} = 1$ and $C$ is an oriented cycle of $D$.
\QED

\begin{Proposition}\label{nwell-div}
Let $f$,$g$ be binomials of $P_D$, such that $f=m_1-m_2$ and $g=m_1-m'_2$. If there is no oriented cycle in $D$ and $m'_2 \mid m_2$ then $f=g$. 
\end{Proposition}
\proof
There exits a monomial $l$ such that $m_2=lm'_2$ since $m'_2 \mid m_2$. On the other hand we have that  $0=\varphi(f)=\varphi(m_1)-\varphi(l)\varphi(m'_2)$, then $\varphi(m_1)=\varphi(l)\varphi(m'_2)$. But $\varphi(m_1)=\varphi(m'_2)$ because $g \in P_D$. Hence, $\varphi(1)=\varphi(l)$ and $(1-l) \in P_D$. Using 
Corollary \ref{oriencycle} and since there is no oriented cycle in $D$ we have that $l=1$ and $f=g$.
\QED

\begin{Corollary}\label{descompcycle}
If $0 \neq t^{\alpha}- t^{\beta}=f$ in $P_D$ with $\gcd(t^{\alpha},t^{\beta})=1$, then there exist cycles $C_1,\dotsc,C_r$ such that $t^{\alpha}=t^{\alpha_1}\dotsm t^{\alpha_s}$ and $t^{\beta}=t^{\beta_1}\dotsm t^{\beta_s}$ where $t_{C_i}=t^{\alpha_i}- t^{\beta_i}$ for $i=1,\dotsc,s$.
\end{Corollary}
\proof
By Theorem \ref{bincycle} there exists a cycle $C_1$ with $t_{C_1}=t^{\alpha_1}- t^{\beta_1}$ such that $t^{\alpha}=t^{\alpha_1}t^{\alpha'}$ and $t^{\beta}=t^{\beta_1}t^{\beta'}$. Since $f=t^{\alpha}- t^{\beta}$ in $P_D$ then $\varphi(f)=0$ and 
\[
\varphi(t^{\alpha_1}) \varphi(t^{\alpha'})=\varphi(t^{\alpha})=\varphi(t^{\beta})=\varphi(t^{\beta_1})\varphi(t^{\beta'}).
\] 
But $t_{C_1} \in P_D$ then $\varphi(t^{\alpha_1})=\varphi(t^{\beta_1})$. Hence $\varphi(t^{\alpha'})=\varphi(t^{\beta'})$ and $f_1=t^{\alpha'} - t^{\beta'} \in P_D$. Since $\gcd(t^{\alpha'},t^{\beta'})=1$, if $f_1=0$ then $t^{\alpha'}=t^{\beta'}=1$ and $f=t_{C_1}$. Now, if $f_1 \neq 0$ by Theorem \ref{bincycle} there exists a cycle $C_2$ with $t_{C_2}=t^{\alpha_2}- t^{\beta_2}$ such that $t^{\alpha'}=t^{\alpha_2}t^{\alpha''}$ and $t^{\beta'}=t^{\beta_2}t^{\beta''}$. Furthermore $t^{\alpha}=t^{\alpha_1}t^{\alpha_2}t^{\alpha''}$, $t^{\beta}=t^{\beta_1}t^{\beta_2}t^{\beta''}$ and $f_2=t^{\alpha''} - t^{\beta''} \in P_D$. We continue this process to obtain the result.
\QED

\begin{Definition}
A binomial $t^{\alpha}- t^{\beta}$ in $P_D$ is called primitive if there exists no other binomial $t^{{\alpha}'}- t^{{\beta}'} \in P_D$ such that $ t^{{\alpha}'}$ divides $ t^{\alpha}$ and $ t^{{\beta}'}$ divides $t^{\beta}$. 
\end{Definition}

\begin{Proposition}
The primitive binomials of $P_D$ are the binomials of $P_D$ associated to cycles. 
\end{Proposition}
\proof
Let $g=t^{\alpha}-t^{\beta}$ be a primitive binomial. By Theorem \ref{bincycle} there exists a cycle $C$ such that $t_C=t^{{\alpha}'}-t^{{\beta}'}$, $t^{\alpha'} \mid t^{\alpha}$ and $t^{\beta'} \mid t^{\beta}$. But $g$ is primitive then $g=t_C$, i.e., $g$ is associated to $C$. Therefore to obtain the result it is only necessary to prove that the binomials associated to cycles are primitive. So, let $t_C=t^{\gamma}-t^{w}$ be the binomial associated to $C$, then gcd($t^{\gamma},t^{w})=1$. Suppose that there exists a binomial $f=t^{\gamma'}-t^{w'} \in P_D$ such that $t^{\gamma'} \mid t^{\gamma}$ and $t^{w'} \mid t^w$. Then, by Theorem \ref{bincycle} there exists a cycle $C'$ such that $t_{C'}=t^{\gamma''}-t^{w''}$, $t^{\gamma''} \mid t^{\gamma'}$ and $t^{w''} \mid t^{w'}$. Hence, $t^{\gamma''} \mid t^{\gamma}$ and $t^{w''} \mid t^ w$, this implies $E(C') \subset E(C)$. Thus, $C'=C$, $t_{C'}=t_{C}$ and $t_C=f$. Therefore $t_C$ is primitive.
\QED

\begin{Proposition}\label{conjgencycl}
$P_D=\left( t_C \mid C \mbox{ is cycle of } G \right)$.
\end{Proposition}
\proof
First we observe that 
\[
t^{\alpha_1}t^{\alpha_2} - t^{\beta_1}t^{\beta_2} = t^{\alpha_1}(t^{\alpha_2}-t^{\beta_2}) + t^{\beta_2}(t^{\alpha_1}-t^{\beta_1})
\]
then $t^{\alpha_1}t^{\alpha_2} - t^{\beta_1}t^{\beta_2} \in (t^{\alpha_1}-t^{\beta_1}, t^{\alpha_2}-t^{\beta_2})$. Now, if $g=t^{\alpha}-t^{\beta} \in P_D$ then by Corollary \ref{descompcycle}, there exist cycles $C_1,\dotsc,C_s$ with $t_{C_i}=t^{\alpha_i}- t^{\beta_i}$ such that $t^{\alpha}=mt^{\alpha'}=mt^{\alpha_1}\dotsm t^{\alpha_s}$ and $t^{\beta}=mt^{\beta'}=mt^{\beta_1}\dotsm t^{\beta_s}$ where $m=\gcd(t^{\alpha},t^{\beta})$. By the first observation we have that $g=m(t^{\alpha'}-t^{\beta'}) \in (t_{C_1},\dotsc,t_{C_s})$. Therefore $P_D=\left( t_C \mid C \mbox{ is cycle of } G \right)$.
\QED

\begin{Lemma}\label{cyclegen}
Let $C_1,C_2$ be two cycles of $D$ whose intersection is an oriented path $P$. Then $C_3=(C_1 \cup C_2)\setminus P^{o}$ is a cycle and $t_{C_3} \in (t_{C_1},t_{C_2})$.
\end{Lemma}
\proof
We can assume that $C_1=(x=x_0,x_1,\dotsc,x_{r_1}=y,y_1,\dotsc,y_{r_2},x)$, $C_2=(x=x_0,x_1,\dotsc,x_{r_1}=y,z_1,\dotsc,z_{r_3},x)$ and $C_1 \cap C_2=P=(x,x_1,\dotsc,x_{r_1-1},y)$ where $(x_i,x_{i+1}) \in E(D)$ for $i=0,\dotsc,r_1-1$, then $C_3=(C_1 \cup C_2)\setminus P^{o} =(y,y_1,\dotsc,y_{r_2},x,z_{r_3},\dotsc,z_1,y)$. Thus, $C_3$ is a cycle. Without loss of generality we can suppose that $t_{i+1}=(x_i,x_{i+1})$ for $i=0,\dotsc,r_1-1$. Hence, $t_{C_1}=t^{\alpha}t^{\alpha_1}-t^{\beta_1}$ and $t_{C_2}=t^{\alpha}t^{\alpha_2}-t^{\beta_2}$ where $\alpha=\sum_{i=1}^{r_1} e_i$. Furthermore $t_{C_3}=t^{\alpha_1}t^{\beta_2}-t^{\alpha_2}t^{\beta_1}$. But $t^{\alpha_1}t^{\beta_2}-t^{\alpha_2}t^{\beta_1}=t^{\alpha_2}t_{C_1}-t^{\alpha_1}t_{C_2}$ therefore, $t_{C_3} \in (t_{C_1},t_{C_2})$.
\QED

\begin{Proposition}\label{jul22-1-05} \rm\cite{Ring} 
$P_D$ is generated by the set of binomials corresponding to cycles without chords.
\end{Proposition}
\proof
By Proposition \ref{conjgencycl} it is only necessary to prove that if 
\[
C=(x_1,x_2,\dotsc,x_s,x_{s+1}=x_1)
\] 
is a cycle of $G$ and $e=\{x_i,x_j\}$ is a chord of $C$ then, $t_C \in (t_{C_1},t_{C_2})$ where $C_1=(x_i,x_{i+1},\dotsc,x_j,x_i)$ and $C_2=(x_1,\dotsc,x_i,x_j,x_{j+1},\dotsc,x_s,x_1)$. But the intersection of $C_1$ and $C_2$ is the edge $e$ then, we obtain the result by Lemma \ref{cyclegen}. 
\QED

\section{$\text{CI}{\cal O}$ graphs and toric ideals}\label{CIO-graphs}

Let $G$ be a graph, with $q$ edges and $r$ connected components such that $V(G)=\{x_1,\dotsc,x_n\}$. For ease of exposition we will say that  $G$ is $\text{CI}{\cal O}$ if the toric ideal $P_{G_{\cal O}}$ is a binomial complete intersection for all orientation ${\cal O}$. Recall that this means that $P_{G_{\cal O}}$ is generated by $q-n+r$ binomials, for every orientation ${\cal O}$ of $G$. The main technical result in this section is Theorem \ref{CIOdelete}, using this theorem and lemma \ref{condelete} we prove that the $\text{CI}{\cal O}$ property is closed under induced subgraphs. Also we prove that a graph is $\text{CI}{\cal O}$ if and only if its blocks are $\text{CI}{\cal O}$. 
\\
\\
Recall that if $v \in V(G)$, the neighborhood of $v$ is the set $N_{G}(v)=\{w \in V(G) \mid \{v,w\} \in E(G)\}$. We take $x \in V(G)$, $G'=G\setminus x$ and ${\cal O}'$ an orientation of $G'$.
If $N_{G}(x)=\{x_1,\dotsc,x_s\}$, then we define the orientation ${\cal O}_x$ of $G$ as follows: set
$E(G_{{\cal O}_x})=E(G'_{{\cal O}'})\cup\{(x,x_1),\dotsc,(x,x_s)\}$.  In the following three lemmas assume that 
$E(G_{{\cal O}_x})=\{t_1,\dotsc,t_q\}$, where $t_i=(x,x_i)$ for $i=1,\dotsc,s$. We denote the oriented graph $G_{{\cal O}_x}$ by $D_x$ (i.e. $D_x=G_{{\cal O}_x}$).

\begin{Lemma}\label{divcycle}
Let $t^\alpha-t^\beta$ be a binomial in $P_{D_x}$. If $\alpha=(\alpha_1,\dotsc,\alpha_q)$ and $\beta=(\beta_1,\dotsc,\beta_q)$, then $\alpha \cdot (\underbrace{1,\dotsc,1}_s,0,\dotsc,0)=\beta \cdot (\underbrace{1,\dotsc,1}_s,0,\dotsc,0)$, i.e., $\sum\limits_{i=1}^s \alpha_i = \sum\limits_{i=1}^s\beta_i$.
\end{Lemma}
\proof
We have $t_i=(x,x_i)$ for $i=1,\dotsc,s$ thus, $\varphi(t_i)=x_i x^{-1}$. Hence, $\varphi(t^{\alpha})=\frac{m_1}{x^{k_1} m_2}$ and $\varphi(t^{\beta})=\frac{n_1}{x^{k_2} n_2}$, where $m_1$,$m_2$, $n_1$ and $n_2$ are monomials. Also set $\gcd(x,m_1)=\gcd(x,m_2)=\gcd(x,n_1)=\gcd(x,n_2)=1$. 
Furthermore $k_1=\sum_{i=1}^s \alpha_i$ and $k_2=\sum_{i=1}^s \beta_i$. On the other hand $t^\alpha - t^\beta \in P_{D_x}$ then $\varphi(t^{\beta})=\varphi(t^{\alpha})$. Therefore, $\frac{m_1}{x^{k_1} m_2} =\frac{n_1}{x^{k_2} n_2}$ and $k_1=\sum_{i=1}^s \alpha_i = \sum_{i=1}^s \beta_i=k_2$.
\QED

\begin{Corollary}\label{div-mon}
Let $t^\alpha-t^\beta$ be a binomial in $P_{D_x}$. If there exists $a \in \{1,\dotsc,s\}$ such that $t_a \mid t^\alpha$ then there exists $b \in \{1,\dotsc,s\}$ such that $t_b \mid t^\beta$.
\end{Corollary}

\begin{Lemma}\label{resvert}
Let $f=t^\alpha-t^\beta$ be a binomial in $P_{D_x} \cap k[t_{s+1},\dotsc,t_q]$. If $B$ is a generating set of binomials of $P_{D_x}$ and $B'=B \cap k[t_{s+1},\dotsc,t_q]$ then $f \in (B')\subset k[t_{s+1},\dotsc,t_q]$. 
\end{Lemma}
\proof
Set $B=\{g_{\lambda}\}_{\lambda \in \Omega}$. Since $P_{D_x}=(B)$, then $f=t^\alpha-t^\beta=\sum_{\lambda \in \Omega}{f_{\lambda} g_{\lambda}}$, where $f_{\lambda} \in k[t_1,\dotsc,t_q]$. We evaluate in the equality above $t_1=t_2=\dotsb=t_s=0$. By Corollary~\ref{div-mon} 
\[
{g_{\lambda}\bigl|_{t_1=t_2=\dotsb=t_s=0}=\left\{\begin{array}{ll}
                                     g_{\lambda} & \mbox{if } g_{\lambda} \in B'\\
                                     0  & \mbox{if } g_{\lambda} \notin B'. 
                                      \end{array}
                                     \right.
}                                         
\]
Hence, we obtain $f=f\bigl|_{t_1=t_2=\dotsb=t_s=0}=\sum_{g_{\lambda} \in B'}{{f'}_{\lambda} g_{\lambda}}$,
where ${f'}_{\lambda}=f_{\lambda} \bigl|_{t_1=t_2=\dotsb=t_s=0} \in k[t_{s+1},\dotsc,t_q]$. Therefore $f\in (B')$.
\QED

\begin{Lemma}\label{Mnempty}
If $G'= G \setminus x$ is a connected graph then for every $j_1,j_2 \in \{1,\dotsc,s\}$ there exists a binomial $t_{j_1}m_{j_1}-t_{j_2}m_{j_2}$ in $P_{D_x}$ with $m_{j_1}, m_{j_2} \in k[t_{s+1},\dotsc,t_q]$.
\end{Lemma}
\proof
We know that $t_{j_1}=(x,x_{j_1})$ and $t_{j_2}=(x,x_{j_2})$. Since $G'$ is connected, there exists a path ${\cal L}_{j_1,j_2}$ in $D'=G'_{{\cal O}'}$ between $x_{j_1}$ and $x_{j_2}$. Thus, $C_{j_1,j_2}=(x,x_{j_1},{\cal L}_{j_1,j_2},x_{j_2},x)$ is a cycle in $D_x$ and $t_{C_{j_1,j_2}}=t_{j_1} m-t_{j_2} m'$ with $m,m'\in k[t_{s+1}, \dotsc, t_q]$. Therefore $t_{C_{j_1,j_2}} \in P_{D_x}$. 
\QED

\begin{Theorem}\label{CIOdelete}
Let $G$ be a connected graph. If $x$ is not a cut vertex of $G$ and $G'=G \setminus x$ is not a $\text{CI}{\cal O}$ graph then $G$ is not a $\text{CI}{\cal O}$ graph.
\end {Theorem}
\proof
Set $\vert V(G)\vert=n$, $\vert E(G)\vert=q$ and $\vert E(G) \setminus E(G') \vert=s$. Since $G'$ is not a $\text{CI}{\cal O}$ graph there exists an oriented graph $D'=G'_{{\cal O}'}$ associated to an orientation ${\cal O}'$ of $G'$ such that $P_{D'}$ is not generated by $\text{ht}(P_{D'})$
binomials. We take the orientation ${\cal O}_x$ of $G$ given by $E(G_{{\cal O}_x})=E(G'_{{\cal O}'})\cup\{(x,x_1),\dotsc,(x,x_s)\}$, where $N_{G}(x)=\{x_1,\dotsc,x_s\}$ and $D_x=G_{{\cal O}_x}$. We can suppose that $t_i=(x,x_i)$ for $1\leq i \leq s$, i.e., $E(D_x)\setminus E(D')=\{t_1,\dots,t_s\}$. We take ${\cal G}=\{f_1,\dotsc,f_k\}$ a minimum generating set of binomials of $P_{D_x}$ and we define ${\cal G'}={\cal G} \cap k[t_{s+1},\dotsc,t_q]$. By Lemma~\ref{resvert}, $P_{D'}=\left( g \mid g \in {\cal G}'\right)$ and $\text{ht}(P_{D'})=(q-s)-(n-1)+1$. Hence $\vert{\cal G}'\vert >(q-s)-(n-1)+1$. Now, we define the auxiliary graph ${\cal H}$ with $V({\cal H})= \{t_1,\dots,t_s\}$ and 
\[
E({\cal H})=\left\{ \{t_i,t_j\} \left\vert  \begin{array}{l}
                              \mbox{There exist monomials } m_i, m_j \in k[t_{s+1},\dotsc,t_q] \mbox{ such that }\\
                              \mbox{    }t_im_i-t_jm_j \in {\cal G} \mbox{ or } t_jm_j-t_im_i \in {\cal G} \mbox{ with } i,j \in \{1,\dotsc,s\} 
                                             
                                             \end{array}
          \right.
          \right\}.
\]
If $A=\left\{t_im_i-t_jm_j \in {\cal G} \mid \{t_i,t_j\}\in E({\cal H}) \right\}$ then $\vert A \vert = \vert E({\cal H}) \vert$. Also $A \cup {\cal G}' \subset {\cal G}$ and $A \cap {\cal G}'=\emptyset$. Hence we have that $\vert E({\cal H})\vert + \vert{\cal G}' \vert\leq \vert{\cal G} \vert $. Now, if ${\cal H}$ is connected then, $\vert E({\cal H}) \vert \geq \vert V({\cal H})\vert-1=s-1$. Furthermore, since $\vert{\cal G}'\vert >(q-s)-(n-1)+1$, then $\vert{\cal G} \vert > (s-1)+(q-s)-(n-1)+1=q-n+1=ht(P_{D_x})$. Thus, $P_{D_x}$ is not generated by $q-n+1$ binomials.
\\
Thus, to obtain the result it is only necessary to prove that ${\cal H}$ is connected. Suppose, by way of contradiction, that the connected components of $\mathcal{H}$ are $\mathcal{H}_1,\dotsc ,\mathcal{H}_p$ with $p>1$. By Lemma \ref{Mnempty} there exists a binomial $t_1m_1-t_jm_j\in P_{D_x}$ such that $t_1\in V({\mathcal H}_1)$ and $t_j \in V({\mathcal H}_2)$, where $m_1,m_j$ are monomials in $k[t_{s+1},\dotsc ,t_q]$ and $j\in\{ 2,\dotsc ,s\}$. Recall that each $f_i\in\mathcal{G}$ has the form $f_i=t^{\alpha_i}-t^{\beta_i}$, where $t^{\alpha_i},t^{\beta_i}$ are monomials in $k[t_1,\dotsc ,t_q]$ then,
\begin{equation*}
t_1m_1-t_jm_j = \sum_{i=1}^k g_i(t^{\alpha_i}-t^{\beta_i}),\qquad\text{with } g_i\in k[t_1,\dotsc ,t_q].
\end{equation*}
If $\alpha_i=(\alpha_1^i,\dotsc ,\alpha_q^i)$ and $\beta_i=(\beta_1^i,\dotsc ,\beta_q^i)$, we take $\alpha_i^\prime=(\alpha_1^i,\dotsc ,\alpha_s^i,0,\dotsc ,0)$ and $\beta_i^\prime=(\beta_1^i,\dotsc ,\beta_s^i,0,\dotsc ,0)$; observe by Lemma~\ref{divcycle}, $\deg(t^{\alpha_i^\prime})=\deg(t^{\beta_i^\prime})$ for $i\in\{ 1,\dotsc ,k\}$. Moreover, $\alpha_i^\prime=0\Leftrightarrow\beta_i^\prime=0$, and in this case, we have that $\bigl(t^{\alpha_i^\prime}-t^{\beta_i^\prime}\bigr)\bigr|_{t_{s+1}=\dotsb =t_q=1} =0$. Hence, if we evaluate in the equation above $t_{s+1}=\dotsb =t_q=1$, we obtain
\begin{equation*}
t_1-t_j = \sum_{\deg(t^{\alpha_i^\prime})=1} g_i^\prime(t^{\alpha_i^\prime}-t^{\beta_i^\prime})+\sum_{\deg(t^{\alpha_i^\prime})>1} g_i^\prime(t^{\alpha_i^\prime}-t^{\beta_i^\prime}),
\end{equation*}
where $g_i^\prime=g_i\bigr|_{t_{s+1}=\dotsb =t_q=1}$. Furthermore, we can write $g_i^\prime=a_i+\lambda_i$ with $a_i\in k[t_1,\dotsc ,t_s]$ and $\lambda_i\in k$, such that if $a_i\ne 0$ then $\deg(a_i)>0$. Thus, 
\begin{align*}
t_1-t_j &= \sum_{\deg(t^{\alpha_i^\prime})=1} \lambda_i(t^{\alpha_i^\prime}-t^{\beta_i^\prime})+g^\prime,
\end{align*}
where $g^\prime = \sum_{\deg(t^{\alpha_i^\prime})=1} a_i(t^{\alpha_i^\prime}-t^{\beta_i^\prime})+\sum_{\deg(t^{\alpha_i})>1} g_i^\prime(t^{\alpha_i^\prime}-t^{\beta_i^\prime})$. Hence 
we have either $\deg(g^\prime)>1$ or $g^\prime=0$. If we suppose $\deg(g^\prime)>1$, then
\begin{align*}
t_1-t_j-\sum_{\deg(t^{\alpha_i^\prime})=1} \lambda_i(t^{\alpha_i^\prime}-t^{\beta_i^\prime}) = g^\prime,
\end{align*} 
however the left side has degree at most one, but $\deg(g^\prime)>1$, which is a contradiction. Therefore $g^\prime=0$, so
\begin{equation*}
t_1-t_j = \hspace{-0.2cm}\sum_{\deg(t^{\alpha_i^\prime})=1} \hspace{-0.2cm}\lambda_i(t^{\alpha_i^\prime}-t^{\beta_i^\prime})= \hspace{-0.2cm}\sum_{\{ t_{i_1},t_{i_2}\}\in E(\mathcal{H}_1)} \hspace{-0.3cm}\lambda_{i_1,i_2}^1(t_{i_1}-t_{i_2})+\dotsb + \hspace{-0.2cm}\sum_{\{ t_{i_1},t_{i_2}\}\in E({\mathcal H}_p)} \hspace{-0.3cm}\lambda_{i_1,i_2}^p(t_{i_1}-t_{i_2}),
\end{equation*}
where $\lambda_i=\lambda_{i_1,i_2}^u$ if $t^{\alpha_i^\prime}=t_{i_1}$, $t^{\beta_i^\prime}=t_{i_2}$ and $\{ t_{i_1},t_{i_2}\}\in E({\mathcal H}_u)$. In the last equation we evaluate $t_\ell=1$ if $t_\ell\in V(\mathcal{H}_1)$ and $t_\ell=0$ otherwise, and we obtain

\begin{equation*}
1-0 = \hspace{-0.3cm}\sum_{\{ t_{i_1},t_{i_2}\}\in E({\mathcal H}_1)} \hspace{-0.3cm}\lambda_{i_1,i_2}^1(1-1)+\hspace{-0.3cm}\sum_{\{ t_{i_1},t_{i_2}\}\in E({\mathcal H}_2)} \hspace{-0.3cm}\lambda_{i_1,i_2}^2(0-0)+\dotsb + \hspace{-0.2cm}\sum_{\{ t_{i_1},t_{i_2}\}\in E ({\mathcal H}_p)} \hspace{-0.3cm}\lambda_{i_1,i_2}^p(0-0)=0
\end{equation*}
since $t_1\in V(\mathcal{H}_1)$ and $t_j \in V(\mathcal{H}_2)$. This is a contradiction. Therefore $\mathcal{H}$ is connected.
\QED

\begin{Definition}
Let $G$ be a graph. A subgraph $H$ of $G$ is a called an induced subgraph of $G$ if $E(H)=\{ e \in E(G) \mid e \subseteq V(H) \}$. Furthermore if $X \subseteq V(G)$ then we denote by $[X]_G$ the induced subgraph $H$ of $G$ such that $V(H)=X$, and in this case we say that $H$ is induced by $X$.
\end{Definition}

\begin{Lemma}\label{condelete}
Let $G$ be a connected graph. If $G'$ is a proper induced connected subgraph of $G$, then there exists $x \in V(G)\setminus V(G')$ such that $G \setminus x$ is connected.
\end{Lemma}
\proof
Set $\vert V(G) \vert =n$. If $v \in V(G) \setminus V(G')$, there exists $G_v$ a connected component of $G \setminus v$ such that $G'\subset G_v$. We define $r_v=\vert V(G_v) \vert$. Let $x$ be the vertex in $V(G)\setminus V(G')$ with $r_x$ maximal. Let $G_1,\dotsc,G_k$ be the connected components of $G \setminus x$ with $G' \subset G_1$, i.e., $G_v=G_1$. Hence, $\vert V(G_1)\vert=r_x$. Suppose that $r_x < n-1$, then $k>1$. Thus, there is $y \in V(G_2) \subset V(G) \setminus V(G')$. Since $G$ is connected, we have $x$ is adjacent to a vertex of $G_1$. Hence $[V(G_1)\cup{x}]_G$ is contained in a connected component of $G\setminus y$. Then $r_y > r_x$, this is a contradiction. Therefore $r_x=n-1$, $k=1$ and $G \setminus x$ is connected.
\QED

\begin{Lemma}\label{CIO-1int}
Let $G_1$, $G_2$ be $\text{CI}{\cal O}$ subgraphs of $G$, such that $G=G_1 \cup G_2$. If $\vert V(G_1) \cap V(G_2) \vert \leq 1$ then $G$ is a $\text{CI}{\cal O}$ graph.
\end{Lemma}
\proof
Let $D=G_{\cal O}$ be the oriented graph associated to an orientation ${\cal O}$ of $G$ and ${\cal O}_i$ the orientation of $G_i$ induced by ${\cal O}$. We set $D_i=(G_i)_{{\cal O}_i}$. Then for $i=1,2$ there exists a binomial generating set ${\cal G}_i$ of $P_{D_i}$ such that $\vert {\cal G}_i \vert=q_i-n_i+r_i$ where $q_i=\vert E(G_i) \vert$, $n_i=\vert V(G_i) \vert$ and $r_i$ is the number of the connected components of $G_i$. Since $\vert V(D_1) \cap V(D_2) \vert \leq 1$, we have that if $C$ is a cycle of $D$ then either $C \subset D_1$ or $C \subset D_2$. Hence, $t_C\in \left({\cal G}_1 \cup {\cal G}_2 \right)$. So ${\cal G}= {\cal G}_1 \cup {\cal G}_2$ is a binomial generating set of $P_{D}$. Since $E(G_1) \cap E(G_2)= \emptyset$ then ${\cal G}_1 \cap {\cal G}_2= \emptyset$ and 
\[
\vert {\cal G} \vert= \vert {\cal G}_1 \vert + \vert {\cal G}_2 \vert = (q_1 + q_2)-(n_1 + n_2)+(r_1 + r_2). 
\]
On the other hand $\vert E(G) \vert = q_1 +q_2$ and the cardinality of $V(G)$ is 
\[
\vert V(G) \vert =\left\{\begin{array}{ll} 
                    n_1 + n_2 & \mbox{if } V(D_1) \cap V(D_2)= \emptyset\\
                    n_1 + n_2 - 1  & \mbox{if } \vert V(D_1) \cap V(D_2) \vert = 1  
             \end{array}
      \right.                                        
\]
but if $r$ is the number of components of $G$ then $r = r_1 + r_2$ in the first case and $r = r_1 + r_2 - 1$
in the second case. Therefore $\vert {\cal G} \vert= \vert E(G) \vert - \vert V(G) \vert + r = \text{ht}(P_D)$   
and $G$ is a $\text{CI}{\cal O}$ graph.
\QED

\begin{Proposition}\label{componentsCIO}
Let $G$ be a graph. $G$ is a $\text{CI}{\cal O}$ graph if and only if its connected components are $\text{CI}{\cal O}$ graphs.
\end{Proposition}
\proof
$\Rightarrow$) By induction on $\vert V(G) \vert$.  Let $G_1, \dotsc, G_k$ be the connected components of $G$. We can suppose that $k > 1$. By Lemma \ref{condelete} there exists $x_i \in V(G_i)$ such that $G_i \setminus x_i$ is connected, for $i=1,\dotsc,k$. Hence by the induction hypothesis the connected components of $G \setminus x_i$ are $\text{CI}{\cal O}$ graphs. Since this is true for every $i=1,\dotsc,k$ then $G_1,\dotsc, G_k$ are $\text{CI}{\cal O}$ graphs. 
\\
$\Leftarrow$) It is a consequence of Lemma \ref{CIO-1int}. 
\QED

\begin{Corollary}\label{inducedCIO}
If $G$ is a $\text{CI}{\cal O}$ graph then every induced subgraph $H$ of $G$ is a $\text{CI}{\cal O}$ graph.
\end{Corollary}
\proof
Assume $H$ is not a $\text{CI}{\cal O}$ graph. By Proposition \ref{componentsCIO}, $H$ has a connected component $H_1$ such that $H_1$ is not a $\text{CI}{\cal O}$ graph. Let $G_1$ be the connected component of $G$ such that $H_1 \subset G_1$. If $G_1 \neq H_1$ then, by Lemma \ref{condelete} there exists $x_1 \in V(G_1)\setminus V(H_1)$ such that $G^1_1=G_1\setminus x_1$ is connected. Furthermore $H_1$ is an induced subgraph of $G^1_1$ thus, if $G^1_1 \neq H_1$ then, there exists $x_2 \in V(G^1_1)\setminus V(H_1)$, such that 
$G^2_1=G^1_1 \setminus x_2$ is a connected subgraph and $H_1 \subset G^2_1$. We continue with this process until we obtain the connected subgraphs $G_1=G^0_1,G^1_1,\dotsc,G^k_1=H_1$ such that $G_1^{i+1}=G_1^i \setminus x_{i+1}$. Since $H_1$ is not a $\text{CI}{\cal O}$ graph by Theorem \ref{CIOdelete}  $H_1=G_1^k,\dotsc,G_1^1,G_1^0=G_1$ are not $\text{CI}{\cal O}$ graphs. Therefore by Proposition \ref{componentsCIO}, $G$ is not a $\text{CI}{\cal O}$ graph. This is a contradiction, therefore $H$ is a $\text{CI}{\cal O}$ graph.
\QED

\begin{Proposition}\label{block-CIO}
$G$ is a $\text{CI}{\cal O}$ graph if and only if every block $B$ of $G$ is a $\text{CI}{\cal O}$ graph.
\end{Proposition}
\proof
$\Rightarrow$) If $B$ is a block of $G$ then $B$ is an induced subgraph of $G$. Hence by Theorem \ref{inducedCIO}, $B$ is a $\text{CI}{\cal O}$ graph.
\\
$\Leftarrow$) It is a consequence of Lemma \ref{CIO-1int}. 
\QED

\section{Chorded-theta subgraphs and transversal triangles}\label{Ch-thetas}

In this section we introduce chorded-theta subgraphs, the notion of transversal triangles in these subgraphs and theta-ring graphs. In Theorem \ref{minimal-thetas} we describe the minimal chorded-theta subgraphs without transversal triangles. In this section we prove that theta-ring graphs are closed under $0,1,2-$clique-sums
and that chordal graphs are theta-ring graphs. With these results we obtain the main result of this section: theta-ring graphs can be constructed by clique-sums of complete graphs and/or cycles, or equivalently by $0,1,2-$clique-sums of chordal graphs and/or cycles. Finally we prove that the minimal forbiden induced subgraphs for the characterization of theta-ring graphs are: prisms, pyramids, thetas and partial wheels that are chorded thetas. This section is independent of the rest of the paper and its techniques are purely combinatorial.

\begin{Definition}
A chorded-theta subgraph $T$ in $G$ is a subgraph induced by three paths ${\cal L}_1$, ${\cal L}_2$, ${\cal L}_3$ each between the non adjacent vertices $x$ and $y$ such that $V({\cal L}_i) \cap V({\cal L}_j)=\{x,y\}$ for $1 \leq i < j \leq 3$. The edges that do not belong to any of the sets $E({\cal L}_1)$, $E({\cal L}_2)$ and $E({\cal L}_3)$ are called the chords of $T$. Furthermore if each chord of $T$ has its end vertices in different paths (of ${\cal L}_1, {\cal L}_2, {\cal L}_3$), then $T$ is called a simple chorded-theta subgraph.
\end{Definition}

\noindent
To emphasize that ${\cal L}_1, {\cal L}_2$ and ${\cal L}_3$ are the three paths associated to $T$ we denote them by ${\cal L}_i (T)$ for $i=1,2,3$. On the other hand we denote by ${\cal L}_i^o$ the interior of ${\cal L}_i$, i.e., ${\cal L}_i^o={\cal L}_i \setminus \{x,y\}$.

\begin{Definition}
Let $T$ be a chorded-theta of $G$, with ${\cal L}_1(T)=(x,x_1,\dotsc,x_{r_1},y)$, ${\cal L}_2(T)=(x,y_1,\dotsc,y_{r_2},y)$ and ${\cal L}_3(T)=(x,z_1,\dotsc,z_{r_3},y)$. A transversal triangle $H$ of $T$ is a triangle in $G$ such that $V(H)=\{x_i,y_j,z_k\}$ for some $i,j,k$.
\end{Definition}

\noindent
We define the following sets $\Theta(G)=\{T \mbox{ }\vert \mbox{ } T \mbox{ is a chorded-theta of } G \}$,
$\Theta_{\Delta}(G)=\{T \in \Theta(G) \mbox{ }\vert \mbox{ } T \mbox{ has a transversal triangle}\}$ and 
$\Theta^c_{\Delta}(G)=\Theta(G) \setminus \Theta_{\Delta}(G)$. Moreover a minimum chorded-theta subgraph $T$ without transversal triangles is a subgraph in $\Theta^c_{\Delta}(G)$ such that if $T' \in \Theta^c_{\Delta}(G)$ then $\vert V(T) \vert \leq \vert V(T') \vert$.

\begin{Remark}\label{minimum-simple}
If $T$ is a minimum chorded-theta subgraph without transversal triangles, then $T$ is a simple chorded-theta subgraph. 
\end{Remark}

\noindent
In the following two lemmas we assume that $T$ is a chorded-theta subgraph of $G$ with ${\cal L}_1(T)=(x,x_1,\dotsc,x_{r_1},y)$, ${\cal L}_2(T)=(x,y_1,\dotsc,y_{r_2},y)$ and ${\cal L}_3(T)=(x,z_1,\dotsc,z_{r_3},y)$, where $\vert V(T) \vert = r_1 + r_2 + r_3 + 2$.

\begin{Lemma}\label{thetalem1}
Let $T$ be a minimum chorded-theta subgraph without transversal triangles with $r_1 \leq r_2$. If $\{x_i,y_j\}$ is a chord of $T$, then $(i,j) \in \{(1,1),(r_1,r_2)\}$ or $r_1=1$.
\end{Lemma}
\proof
By Remark \ref{minimum-simple}, $T$ is a simple chorded-theta. Let $\{x_i,y_j\}$ be a chord of $T$. We can assume that $r_1>1$. Hence, $i<r_1$ or $1<i$.
\\
Now, suppose $i<r_1$ and take the chorded-theta $T_1$ with ${\cal L}_1(T_1)=(x_i,\dotsc,x_{r_1},y)$, ${\cal L}_2(T_1)=(x_i,y_j,\dotsc,y_{r_2},y)$ and ${\cal L}_3(T_1)=(x_i,x_{i-1},\dotsc,x_1,x,z_1,\dotsc,z_{r_3},y)$.
Since $T$ has not chords whose both its end vertices are in ${\cal L}_1(T)$ and $T \in \Theta^c_{\Delta}(G)$, then $T_1 \in \Theta^c_{\Delta}(G)$. Observe that $\vert V(T_1)\vert=\vert V(T) \vert - (j-1)$, so
by the minimality of $\vert V(T) \vert$, we have that $j=1$. We will prove that $i=1$, by contradiction we will assume $i>1$. Thus, we can take the chorded-theta $T_2$ with ${\cal L}_1(T_2)=(x,y_1,x_i)$, ${\cal L}_2(T_2)=(x,x_1,\dotsc,x_i)$ and ${\cal L}_3(T_2)=(x,z_1,\dotsc,z_{r_3},y,x_{r_1},\dotsc,x_i)$. Again, since $T$ has not chords whose both end vertices are in ${\cal L}_1(T)$ and $T \in \Theta^c_{\Delta}(G)$, we have that $T_2 \in \Theta^c_{\Delta}(G)$. But $\vert V(T_2) \vert=\vert V(T) \vert - (r_2-1)$. Furthermore $1<r_1 \leq r_2$ so $\vert V(T) \vert > \vert V(T_2) \vert$ which is a contradiction. Therefore $i=1$ and $(i,j)=(1,1)$.
\\
Similarly, if $1<i$ we obtain that $(i,j)=(r_1,r_2)$. 
\QED

\begin{Lemma}\label{thetalem}
Let $T$ be a minimum chorded-theta without transversal triangles such that $r_1=1$ and $r_2 \leq r_3$. If there exist chords $e_1, e_2 $ of $T$ with $e_1=\{x_1,z_{s_1}\}$ and $e_2=\{y_j,z_{s_2}\}$ then $\{x_1,y_j\} \notin E(G)$. Moreover $(s_1,j,s_2)=(r_3,1,1)$ or $(s_1,j,s_2)=(1,r_2,r_3)$.
\end{Lemma}
\proof
By Remark \ref{minimum-simple}, $T$ is a simple chorded-theta.
\\
Case (1) $s_1 \leq s_2$. Suppose that $\{x_1,y_j\} \in E(G)$. Since $T \in \Theta^c_{\Delta}(G)$ then $\{x_1,z_{s_2}\} \notin E(G)$ and $z_{s_1} \neq z_{s_2}$. Thus, we take a chorded-theta $T'$ given by 
${\cal L}_1(T')=(z_{s_2},y_j,x_1)$, ${\cal L}_2(T')=(z_{s_2},\dotsc,z_{r_3},y,x_1)$ and ${\cal L}_3(T')=(z_{s_2},z_{s_2-1},\dotsc,z_{s_1},x_1)$. Furthermore there is no chord whose both end vertices are in ${\cal L}_3(T)$, hence $T' \in \Theta^c_{\Delta}(G)$. But $\vert V(T') \vert= \vert V(T)\vert -s_1-r_2+1 < \vert V(T) \vert$ which is a contradiction. Therefore $\{x_1,y_j\} \notin E(G)$. Now, we consider the chorded-theta $T''$ given by ${\cal L}_1(T'')=(x_1,z_{s_1},z_{s_1+1},\dotsc,z_{s_2},y_j)$, ${\cal L}_2(T'')=(x_1,x,y_1,\dotsc,y_j)$ and ${\cal L}_3(T'')=(x_1,y,y_{r_2},y_{r_2-1},\dotsc,y_j)$. We have that $T'' \in \Theta^c_{\Delta}(G)$ because there is no chord with its both end vertices in ${\cal L}_2(T)$. But 
$\vert V(T'') \vert= \vert V(T) \vert - (s_1 - 1) - (r_3-s_2)$ then, $s_1=1$ and $r_3=s_2$ because $T$ is minimum. If $j \neq r_2$ we can take a chorded-theta subgraph $T'''$, where ${\cal L}_1(T''')=(y_j,\dotsc,y_{r_2},y)$, ${\cal L}_2(T''')=(y_j,\dotsc,y_1,x,x_1,y)$ and ${\cal L}_3(T''')=(y_j,z_{r_3},y)$. Thus $T''' \in \Theta^c_{\Delta}(G)$ because $T$ is a simple chorded-theta together with the fact 
that $T \in \Theta^c_{\Delta}(G)$. But $\vert V(T''')\vert=\vert V(T) \vert - (r_3 - 1)$ hence $r_3=1$. Since $r_2 \leq r_3$ then $r_2=1=j$, which is a contradiction. Therefore $j=r_2$ and $(s_1,j,s_2)=(1,r_2,r_3)$.
\\
Case (2) $s_2 \leq s_1$. Similarly to the case (1). We obtain that $\{x_1,y_j\} \notin E(G)$ and $(s_1,j,s_2)=(r_3,1,1)$. 
\QED

\begin{Definition}\rm{[$\forall \theta \exists \Delta-$property]}
A graph $G$ is called a theta-ring graph if every chorded-theta of $G$ has a transversal triangle, in this
case we say that $G$ has the $\forall \theta \exists \Delta-$property. 
\end{Definition}

\begin{Remark}
$G$ is a theta-ring graph if and only if $\Theta^c_{\Delta}(G)=\emptyset$.
\end{Remark}

\begin{Theorem}\label{minimal-thetas}
Let $G$ be a graph. If $G$ is not a theta-ring graph then there exists a chorded-theta $T$ in $G$ where ${\cal L}_1(T)=(x,x_1,\dotsc,x_{r_1},y)$, ${\cal L}_2(T)=(x,y_1,\dotsc,y_{r_2},y)$, ${\cal L}_3(T)=(x,z_1,\dotsc,z_{r_3},y)$, with $1 \leq r_1 \leq r_2 \leq r_3$ such that $T$ satisfies at least one of the following conditions:
\begin{description}
\item{\rm (1)} The chords of $T$ are contained in  $\left \{ \{y_1,z_1\},\{y_{r_2},z_{r_3}\} \right \}$.
\item{\rm (2)} $r_1=1$, and the chords of $T$ are contained in $\left \{ \{x_1,z_{r_3}\},\{y_1,z_1\} \right \}$.
\item{\rm (3)} $r_1=1$, and the chords of $T$ have the form $\{x_1,z_j\}$ for some $j$.
\item{\rm (4)} $r_1=1$, $r_2=1$ and the chords of $T$ have the form $\{x_1,y_1\}$ or $\{x_1,z_j\}$
for some j.
\end{description}
\end{Theorem}
\proof
Since $\Theta^c_{\Delta}(G)\neq \emptyset$. We can take a minimum chorded-theta $T$ without transversal triangles. Thus, $T$ is a simple chorded-theta. Let ${\cal L}_1(T)=(x,x_1,\dotsc,x_{r_1},y)$, ${\cal L}_2(T)=(x,y_1,\dotsc,y_{r_2},y)$ and ${\cal L}_3(T)=(x,z_1,\dotsc,z_{r_3},y)$ with $1 \leq r_1 \leq r_2 \leq r_3$. We can suppose that $\vert{\cal L}_1(T)\vert \leq \vert{\cal L}_1(T')\vert$ for all $T' \in \Theta^c_{\Delta}(G)$ such that $\vert V(T') \vert$ is minimum. Also we assume that if $\vert {\cal L}_1(T)\vert=\vert {\cal L}_1(T') \vert $ then $\vert{\cal L}_2(T)\vert \leq \vert{\cal L}_2(T')\vert$. 
\\
Now, we suppose that $T$ does not satisfy property \rm (1). Hence, there exists a chord of $T$ with one of the possible forms: $\{x_l,y_j\}$ or $\{x_l,z_j\}$ or $\{y_i,z_k\}$, where $(i,k)\notin \{(1,1),(r_2,r_3)\}$. 
\\
First case, we suppose that there exists $e=\{y_i,z_k\}$ a chord of $T$, where $(i,k)\notin \{(1,1),(r_2,r_3)\}$. By Lemma \ref{thetalem1} we have that  $r_2=1$ and $i=1$. Now, since $r_1 \leq r_2$ then $r_1=1$. So, by Lemma \ref{thetalem}, if there exists $\{x_1,z_j\} \in E(G)$ then $\{x_1,y_1\} \notin E(G)$ and $(j,i,k)=(r_3,1,1)$ or $(j,i,k)=(1,1,r_3)$. Thus, we would be in case (2). Now, if $\{x_1,z_j\} \notin E(G)$ for every $j$ then the chords of $T$ have the form  $\{x_1,y_1\}$ or $\{y_1,z_j\}$, and we are in case (4).
\\
Second case, suppose there exists a chord $e$ of $T$ such that $e=\{x_i,y_j\}$ or $e=\{x_i,z_j\}$. We will prove that $r_1=1$. We argue by contradiction, assume that $r_1>1$ and $e=\{x_i,y_j\}$. By Lemma \ref{thetalem1}, we have that $(i,j) \in \{(1,1),(r_1,r_2)\}$. We can suppose that $(i,j)=(r_1,r_2)$ and take the chorded-theta graph $T_1$ where ${\cal L}_1(T_1)=(x,x_1,\dotsc,x_{r_1})$, ${\cal L}_2(T_1)=(x,y_1,\dotsc,y_{r_2},x_{r_1})$ and ${\cal L}_3(T_1)=(x,z_1,\dotsc,z_{r_3},y,x_{r_1})$. 
Since $T$ is a simple chorded-theta and $T \in \Theta^c_{\Delta}(G)$ then $T_1 \in \Theta^c_{\Delta}(G)$. But $\vert V(T) \vert=\vert V(T_1) \vert$ and $\vert {\cal L}_1(T_1) \vert = \vert {\cal L}_1(T) \vert -1 $, which is a contradiction. The same result is obtained if $e=\{x_i,z_j\}$. Hence, we obtain that $r_1=1$. 
\\
Now, we will prove that if there exists a chord $e'$ of $T$ such that $e'=\{x_1,y_u\}$ then $T$ satisfies (4). First we will prove that $r_2=1$. Again we argue by contradiction, assume $r_2>1$. Hence, $u>1$ or $u<r_2$. Without loss of generality, take $u>1$ and consider the chorded-theta $T_2$ where ${\cal L}_1(T_2)=(x,x_1,y_u)$, ${\cal L}_2(T_2)=(x,y_1,\dotsc,y_u)$ and ${\cal L}_3(T_2)=(x,z_1,\dotsc,z_{r_3},y,y_{r_2},\dotsc,y_u)$. Since $T \in \Theta^c_{\Delta}(G)$ and there are not chords in ${\cal L}_2 (T_2)$, then $T_2 \in \Theta^c_{\Delta}(G)$. But $\vert V(T) \vert = \vert V(T_2) \vert$, $\vert {\cal L}_1(T) \vert = \vert {\cal L}_1(T_2) \vert$ and $\vert {\cal L}_2(T) \vert > \vert {\cal L}_2(T_2) \vert$, which is a contradiction. Thus, $r_1=r_2=1$ and $e'=\{x_1,y_1\}$. Suppose $T$ does not satisfy (4) then there must exist $e_1,e_2 \in E(G)$ such that $e_1=\{x_1,z_{j_1}\}$ and $e_2=\{y_1,z_{j_2}\}$. But, by Lemma \ref{thetalem}, $\{x_1,y_1\} \notin E(G)$. This is a contradiction. Hence, $T$ satisfies (4).
\\
Now, we can  assume that $T$ does not have a chord of the form $e'=\{x_1,y_u\}$ then $e=\{x_1,z_j\}$. Suppose $T$ does not satisfy (3), then there must exists $e_3=\{y_i,z_k\}$. By Lemma \ref{thetalem}, $(j,i,k) \in \{(r_3,1,1),(1,r_2,r_3)\}$. Therefore $T$ satisfies (2).
\QED

\begin{Definition}
A partial wheel $W$ is a graph where $V(W)=\{z,z_1,\dotsc,z_k\}$ such that $C=(z_1,\dotsc,z_k,z_1)$ is a cycle in $W$ and the edges of $W$ are the edges of $C$ and some edges between $z$ and vertices of $C$. $C$ is called the rim of $W$ and $z$ is called the center of $W$.
\\
A prism is a graph consisting of two vertex-disjoint triangles $C_1=(x_1,x_2,x_3,x_1)$ and $C_2=(y_1,y_2,y_3,y_1)$, and three paths $L_1,L_2,L_3$ pairwise vertex-disjoint, such that each $L_i$ is a path between $x_i$ and $y_i$ for $i=1,2,3$ and the subgraph induced by  $V(L_i) \cup V(L_j)$ is a cycle for $1 \leq i < j \leq 3$.
\\
A pyramid is a graph consisting of a vertex $w$, a triangle $C=(z_1,z_2,z_3,z_1)$, and three paths $P_1,P_2,P_3$, such that: $P_i$ is between $w$ and $z_i$ for $i=1,2,3$; $V(P_i) \cap V(P_j)=\{w\}$ and the subgraph induced by $V(P_i) \cup V(P_j)$ is a cycle for $1 \leq i < j \leq 3$; and at most one of the $P_1,P_2,P_3$ has only one edge.
\\
A theta is a graph consisting of two non adjacent vertices $x$ and $y$, and three paths $P,Q,R$ with ends $x$ and $y$, such that the union of every two of $P,Q,R$ is an induced cycle.
\end{Definition} 

\noindent
Detecting prisms, pyramids and thetas has been widely studied, for example see \cite{C-K,LLMT} and their references.

\begin{Definition}
A partial wheel $T$ with rim $C$ and center $z$ is called a $\theta-$partial wheel if $|V(C)|\geq 4$ and there exist two non adjacent vertices in $V(C)\cap N_T(z)$.
\end{Definition}

\begin{Remark}
A theta is a chorded-theta without chords. Furthermore, if $W$ is a partial wheel then $W$ is a $\theta-$partial wheel if and only if $W$ is a chorded-theta.
\end{Remark}

\begin{Corollary}\label{min-fobb-thetas}
$T$ is a minimal forbidden induced subgraph for the class of theta-ring graphs if and only if $T$ satisfies at least one of the following conditions:
\begin{description}
\item{\rm a)} $T$ is a $\theta-$partial wheel.
\item{\rm b)} $T$ is a prism, pyramid or theta.
\end{description}
\end{Corollary}
\proof
First we will prove that if $T$ satisfies a) or b) then, $T$ is a chorded-theta graph without transversal triangles and $T$ does not contain a proper chorded-theta subgraph.
\\
If $T$ satisfies a) with rim $C=(z_1,z_2,\dotsc,z_{r_1},z_1)$ and center $z$ then, there exist $z_i,z_j \in N_T(z)$ where $\{z_i,z_j\} \notin E(T)$. We can suppose that $i <j$ then, $T$ is a chorded-theta graph with the paths ${\cal L}_1(T)=(z_i,z_{i+1},\dotsc,z_j)$, ${\cal L}_2(T)=(z_i,z,z_j)$ and ${\cal L}_3(T)=(z_i,z_{i-1},\dotsc,z_1,z_{r_1},z_{{r_1}-1},\dotsc,z_j)$. Furthermore, since $C$ does not contain chords then $T$ is a chorded-theta without transversal triangles. On the other hand if $T'$ is a chorded-theta subgraph of $T$, then one terminal vertex is in $C$. Thus, one path of $T'$ has the form $(z_{i_1},z,z_{j_1})$ with $\{z_{i_1},z_{j_1}\} \notin E(T)$ and $i_1 < j_1$. Hence, the others paths of $T'$ are $(z_{i_1},z_{{i_1}+1},\dotsc,z_{j_1})$ and $(z_{i_1},z_{{i_1}-1},\dotsc,z_1,z_{r_1},z_{{r_1}-1},\dotsc,z_{j_1})$. Therefore $V(T)=V(T')$ and $T=T'$.
\\
Now, if $T$ is a prism then, $T$ is a chorded-theta with ${\cal L}_1(T)=(x_1,x_2,L_2,y_2)$, ${\cal L}_2(T)=(x_1,L_1,y_1,y_2)$ and ${\cal L}_3(T)=(x_1,x_3,L_3,y_3,y_2)$. Furthermore $T$ has only two chords, then $T$ does not contain transversal triangles. On the other hand, if $T'$ is a chorded-theta in $T$ then, the terminal vertices are $x_i$ and $y_j$ for some $i,j \in \{1,2,3\}$. Since $x_i$ and $y_j$ have degree at least three in $T'$, we have that $V(T')=V(T)$ and $T'=T$.
\\
If $T$ is a pyramid then, we can suppose that $\vert E(P_1) \vert \geq 2$ and $T$ is a chorded-theta with  
${\cal L}_1(T)=(z_1,P_1,w)$, ${\cal L}_2(T)=(z_1,z_2,P_2,w)$ and ${\cal L}_3(T)=(z_1,z_3,P_3,w)$. Additionally $T$ has only one chord then $T$ does not contain transversal triangles. Furthermore if $T'$ is a chorded-theta in $T$ then the terminal vertices are $z_i$ and $w$ for some $i \in \{1,2,3\}$. Since $z_i$ and $w$ have degree at least $3$ in $T'$, we have that $V(T')=V(T)$ and $T'=T$.
\\
Finally, if $T$ is a theta graph it is clear that $T$ is a minimal chorded-theta graph without transversal triangles.
\\
Now, to obtain the result is only necessary to prove that if $T$ satisfies (1), (2), (3) or (4) of Theorem \ref{minimal-thetas} then $T$ satisfies a) or b).
\\
Case $T$ satisfies (1). If $T$ has less than two chords and $T$ is neither a pyramid nor a theta, then
$r_3=1$ and $r_2=r_1=1$. Thus, $T$ satisfies a). Now, we can suppose that $T$ has exactly two chords. If $r_2 \geq 2$ then $T$ is a prism. Furthermore if $r_2=1$ then $T$ satisfies a).
\\
Case $T$ satisfies (2). If $T$ has less than two chords, and $T$ is neither a pyramid nor a theta then $r_3=1$ and $r_2=1$. Hence, $T$ satisfies a). Now, we can suppose that $T$ has exactly two chords. If $r_3 \geq 2$ then $T$ is a prism. Furthermore if $r_3=1$ then $r_2=r_1=1$ and $T$ satisfies a). 
\\
If $T$ satisfies (3) or (4) then $T$ satisfies a).
\QED

\begin{Definition}
Let $G_1,G_2$ be graphs such that $G=G_1 \cup G_2$, $K=G_1 \cap G_2$. If $K$ is a complete graph with $\vert  K \vert=k$, then $G$ is called the $k-$clique-sum (or clique-sum) of $G_1$ and $G_2$ in $K$.
\end{Definition}

\begin{Remark}
Let $r_i$ be the number of connected components of $G_i$ for $i=1,2$. If $G$ is $0-$clique-sum of $G_1$ and $G_2$ then the number of connected components of $G$ is $r_1+r_2$.
\end{Remark}

\noindent
In our definition of clique-sum we do not allow to delete any of the edges of $K$.

\begin{Proposition}\label{sum-thetaring}
The clique-sum of theta-ring graphs is a theta-ring graph.
\end{Proposition}
\proof
Let $G_1$ and $G_2$ be theta-ring graphs. We assume that $G$ is a $k-$clique-sum of $G_1$ and $G_2$. We can suppose that $G=G_1 \cup G_2$ (i.e., $G_1,G_2 \subseteq G$) and $G_1 \cap G_2=K$ where $K$ is a $k-$complete subgraph of $G$. By contradiction suppose that $G$ is not a theta-ring graph. Hence, there exists a chorded-theta $H$ (of $G$) without transversal triangles, where ${\cal L}_1(H)=(x,x_1,\dotsc,x_{r_1},y)$, ${\cal L}_2(H)=(x,y_1,\dotsc,y_{r_2},y)$ and ${\cal L}_3(H)=(x,z_1,\dotsc,z_{r_3},y)$. We can assume that $H$ is a simple chorded-theta. Since $\{x,y\} \notin E(G)$ then $\{x,y\} \cap (G_1 \triangle G_2) \neq \emptyset$, where $G_1 \triangle G_2$ is the symmetric difference. Thus, without loss of generality we can assume that $x \in V(G_1) \setminus V(G_2)$. If $y \in V(G_2) \setminus V(G_1)$ then ${\cal L}^{o}_i(H) \cap K \neq \emptyset$ for $i=1,2,3$. But $K$ is a complete graph then $H$ has a transversal triangle, which is not possible so $y \in V(G_1)$. On the other hand $G_1$ is a theta-ring graph then, $V(H)$ is not contained in $V(G_1)$. Hence, there exists $w \in V(H)$ such that $w \in V(G_2) \setminus V(G_1)$. We can suppose that $w \in {\cal L}^{o}_1(H)$, i.e., $w=x_i$ for some $1\leq i \leq r_1$. Since $H$ is a simple chorded-theta then, $L_1=(x,x_1,\dotsc,x_i=w)$
is a path between $G_1 \setminus G_2$ and $G_2 \setminus G_1$. Furthermore $K$ is a cutset then there exists $x_{i_1} \in V(L_1) \cap V(K)$ with $i_1 < i$. On the other hand if $y \in V(K)$ then $\{y,x_{i_1}\} \in E(G)$ but this is not possible because $H$ is a simple chorded-theta. Thus, $y \in V(G_1) \setminus V(G_2)$ and $L_2=(w=x_i,x_{i+1},\dotsc,x_{r_1},y)$ is a path between $G_1 \setminus G_2$ and $G_2 \setminus G_1$. Hence, there exists $x_{i_2} \in V(L_2) \cap V(K)$ with $i < i_2$. Since, $x_{i_1}, x_{i_2} \in V(K)$ then $\{x_{i_1},x_{i_2}\} \in E(G)$. This is a contradiction, because $H$ is a simple chordad-theta. 
\QED

\begin{Lemma}\label{Paths-connected1}
Let $G$ be a $2-$connected graph and let $H$ be an induced connected subgraph of $G$ with $\vert V(H) \vert \geq 2$. If $x \in V(G) \setminus V(H)$ then there exist paths $L_1$ and $L_2$ between $x$ and $H$ such that $V(L_1) \cap V(L_2)=\{x\}$, $V(L_i)\cap V(H)=\{a_i\}$ for $i=1,2$, and $a_1 \neq a_2$.
\end{Lemma}
\proof 
Let $y$ be a vertex in $H$. Since $G$ is $2-$connected graph there exist paths $P_1=(x=x_1,x_2,\dotsc,x_{s_1}=y)$ and $P_2=(x=y_1,y_2,\dotsc,y_{s_2}=y)$ between $x$ and $y$, such that $V(P_1) \cap V(P_2)=\{x,y\}$. We define $j_1=$ min$\{k \mid x_k \in V(H)\}$ and $j_2=$ min$\{k \mid y_k \in V(H)\}$. If $x_{j_1} \neq y_{j_2}$ we can take $L_1=(x=x_1,x_2,\dotsc,x_{j_1})$ and $L_2=(x=y_1,y_2,\dotsc,y_{j_2})$. Hence, we can assume that $x_{j_1}=y_{j_2}$, but $V(P_1) \cap V(P_2) = \{x,y\}$ then $x_{j_1}=y_{j_2}=y$.
Since $\vert V(H) \vert \geq 2$ thus there must exist $a \in V(H) \setminus \{y\}$. Furthermore $G$ is $2-$connected so $G\setminus y$ is connected, and there exists a path $P_3=(x=z_1,\dotsc,z_{s_3}=a)$ between $x$ and $a$ in $G \setminus y$. Now, we take $k_1=$ max$\{k \mid z_k \in V(P_1) \cup V(P_2)\}$, $k_2=$ min$\{k \mid k \geq k_1 \mbox{ and } z_k\in V(H)\}$. We can suppose that $z_{k_1} \in V(P_1)$ and $z_{k_1}=x_{i_1}$. Finally we take $L_1=(x=x_1, \dotsc, x_{i_1}=z_{k_1},z_{{k_1}+1},\dots,z_{k_2})$ and $L_2=P_2$.
\QED

\begin{Definition}
Let $C$ be a cycle without chords. $C$ is called a hole if $C$ is not a triangle.
\end{Definition}

\begin{Lemma}\label{Paths-connected2}
Let $G$ be a theta-ring graph with a hole $C$. Let $x$ be a vertex of $G$ such that $x \notin V(C)$.
If there exist paths $P_1$ and $P_2$  between $x$ and $C$ such that $V(P_1)\cap V(P_2)=\{x\}$ with $V(P_i)\cap V(C)=\{y_i\}$ for $i=1,2$ then, $\{y_1,y_2\} \in E(C)$. Furthermore if $P$ is a path between $x$ and $C$ with $V(P)\cap V(C)=\{y\}$ then $y \in \{y_1,y_2\}$.
\end{Lemma}
\proof
Set $C=(x_1,x_2,\dotsc,x_k=x_1)$ where $x_1=y_1$ and $y_2=x_s$. By way of contradiction, suppose that $\{y_1,y_2\} \notin E(C)$. Thus, there exists a chorded-theta $H$ of $G$ with ${\cal L}_1(H)=(y_1,P_1,x,P_2,y_2)$, ${\cal L}_2(H)=(y_1=x_1,x_2,\dotsc,x_s=y_2)$ and ${\cal L}_3(H)=(y_1=x_1,x_k,x_{k-1},\dotsc,x_s=y_2)$ but $C$ does not have any chord then, $H$ has no transversal triangles. Hence, $\{y_1,y_2\} \in E(C)$. Now, if $P=(y=z_1,z_2,\dotsc,z_{s_2}=x)$ is a path between $x$ and $C$ with $V(P)\cap V(C)=\{y\}$, we will prove that $y \in \{y_1,y_2\}$. By way of contradiction suppose that $y \notin \{y_1,y_2\}$. We take $k_3=$ min$\{i \mid z_i \in V(P_1) \cup V(P_2)\}$, but $y \notin V(P_1) \cup V(P_2)$ so $k_3 > 1$. We can suppose that $z_{k_3} \in V(P_1)$. Set $P_1=(w_1=y_1,w_2,\dotsc,w_{s_1}=x)$ and
$z_{k_3}=w_{l}$. Thus, there exist paths $L_1=(y_1=w_1,w_2,\dotsc, w_l=z_{k_3})$, $L_2=(y=z_1,z_2, \dotsc,z_{k_3})$ and $L_3=(y_2,P_2,x=w_{s_1},w_{s_1-1},\dotsc,w_l=z_{k_3})$ such that $V(L_1) \cap V(C)=\{y_1\}$, $V(L_2) \cap V(C)=\{y\}$ and $V(L_3) \cap V(C)=\{y_2\}$. Then, $\{y,y_1\},\{y,y_2\} \in E(C)$, but $\{y_1,y_2\} \in E(C)$ and $C$ is a hole, which is a contradiction. Therefore, $y \in \{y_1,y_2\}$.
\QED

\begin{Definition}
A graph is chordal if each of its cycles of length at least $4$ has a chord.
\end{Definition}

\begin{Proposition}\label{sumcomplchd}\cite{Seymour-Weaver}
Chordal graphs are exactly the graphs that can be formed by clique-sums of complete graphs.
\end{Proposition}

\begin{Theorem}\label{For-all-Theta-exists-delta}\rm{[$\forall \theta \exists \Delta$]}
Let $G$ be a graph. The following conditions are equivalent:
\begin{itemize}
\item[i)] G is a theta-ring graph (i.e., $\Theta^c_{\Delta}(G)=\emptyset$).
\item[ii)] $G$ can be constructed by $0,1,2-$clique-sums of chordal graphs and/or cycles. 
\item[iii)] $G$ can be constructed by clique-sums of complete graphs and/or cycles.
\item[iv)] $G$ does not contain as induced subgraph any graph from the following families:
\begin{itemize}
\item[a)] $\theta-$partial wheels.
\item[b)] Prisms, pyramids and thetas.
\end{itemize}
\end{itemize}
\end{Theorem}
\proof
i) $\Leftrightarrow$ iv) By Corollary \ref{min-fobb-thetas}.
\\
iii) $\Rightarrow$) i) Complete graphs and cycles are theta-ring graphs (since, they do not contain chorded-thetas). Then we finish by Proposition \ref{sum-thetaring}.
\\
\\
ii) $\Rightarrow$) iii) Suppose that $G$ can be constructed by $0,1,2-$clique-sums of $H_1, \dotsc, H_s$, 
where each $H_i$ is a chordal graph or cycle. If $s=1$ then, $G$ is a cycle or chordal graph. Hence, $G$ satisfies iii) by Proposition \ref{sumcomplchd}. Now, we can suppose that $s>1$ and we take $G'$ the
induced subgraph obtained by $0,1,2-$clique-sum of $H_1, \dotsc, H_{s-1}$ such that $G$ is the $j-$clique-sum of $G'$ and $H_s$ where $j \in \{0,1,2\}$. Furthermore $\vert V(G') \cap V(H_s) \vert \leq j$. By induction hypothesis we have that $G'$ can be constructed by clique-sums of $L_1, \dotsc, L_{l_1}$, where each $L_i$ is a cycle or a complete subgraph of $G'$. If $H_s$ is a cycle then we obtain the result. Hence, we
can assume that $H_s$ is a chordal graph. Thus, by Proposition \ref{sumcomplchd}, we have that $H_s$ can be constructed by clique-sums of the complete subgraphs $K_1,K_2,\dotsc,K_{l_2}$. Since $\vert V(G') \cap V(H_s) \vert \leq j$ we can suppose that $G' \cap H_s \subset K_u$ for some $u \in \{1, \dotsc, l_2\}$. Therefore $G$ can be constructed by clique-sums of $L_1,\dotsc,L_{l_1},K_{u},K_{u+1},\dotsc,K_{l_2},K_{u-1},\dotsc,K_1$.
\\
\\
i) $\Rightarrow$) ii) By induction on $\vert V(G) \vert$. We can suppose that $\vert V(G) \vert > 1$. Let $B_1,\dotsc,B_s$ be the blocks of $G$. If $s>1$ then by induction hypothesis $B_1,\dotsc,B_s$ satisfy ii). But $G$ can be obtained as $0$ or $1-$clique-sums of its blocks so $G$ satisfies ii). Thus, we can assume that $s=1$ and $G=B_1$. We can suppose that $G$ is not a chordal graph. Thus, $G$ has a cycle $C$ without chords such that $C$ is not a triangle, i.e., $C$ is a hole. Since $G$ is a block we have that $G$ is $2-$connected. Set $C=(x_1,\dotsc,x_k,x_{k+1}=x_1)$ and we define 
\[
C_i=
   \left\{ x\in V(G\setminus C) \left| 
      \begin{array}{l}
         \mbox{There exist paths $L_x, L_x^\prime$ between $x$ and $C$ such that}\\
         \mbox{$V(L_x)\cap V(L_x^\prime)=\{ x\}$, $V(L_x)\cap V(C)=\{ x_i\}$}\\ 
         \mbox{and $V(L_x^\prime)\cap V(C)=\{ x_{i+1}\}$}
      \end{array}
   \right.
   \right\}    
\]
for $i=1,\dotsc,k$. We will prove that $C_1,\dotsc,C_k$ is a partition of $V(G \setminus C)$. By Lemma \ref{Paths-connected1} if $x\in V(G \setminus C)$ then there exist paths $L_1$ and $L_2$ between $x$ and $C$, such that $V(L_1)\cap V(L_2)=\{x\}$ and $V(L_i)\cap V(C)=\{y_i\}$ for $i=1,2$ with $y_1 \neq y_2$. Hence, by Lemma \ref{Paths-connected2}, $\{y_1,y_2\} \in E(C)$. Thus, there exists $i \in \{1,\dotsc,k\}$ such that $x_i=y_1$ and $y_2=x_{i-1}$ or $y_2=x_{i+1}$. So, $x \in V(C_i)$ or $x \in V(C_{i-1})$. Therefore $V(G \setminus C) = \cup^{n}_{i=1} V(C_i)$. Now, we take $z \in V(C_i) \cap V(C_j)$ with $i \neq j$. By Lemma \ref{Paths-connected2}, we have that $\{x_i,x_{i+1}\} \subset \{x_j,x_{j+1}\}$ and $\{x_j,x_{j+1}\} \subset \{x_i,x_{i+1}\}$. Then, $i=j$ which is a contradiction. Hence, $V(C_i) \cap V(C_j) = \emptyset$ and $V(C_1), \dotsc, V(C_r)$ is a partition of $V(G \setminus C)$. Let $G_i$ be the subgraph induced by $V(C_i) \cup \{x_i,x_{i+1}\}$. By the definition of $C_i$ we have that $V(L_x)\cup V(L'_x) \subset V(G_i)$ for $x \in V(C_i)$, so $G_i$ is $2-$connected, and by the induction hypothesis $G_i$ can be obtained as 0,1,2-clique-sums of chordal graphs and cycles. We define $B_1$ as the $2-$clique-sum of $C$ and $G_1$ at the edge $\{x_1,x_2\}$.
Now, we will construct successively $B_{i+1}$ as the $2-$clique-sum of $B_i$ and $G_{i+1}$ at the edge $\{x_{i+1},x_{i+2}\}$ for $i=1,\dotsc,k-2$. But $B_{k-1}=G$ therefore $G$ can be constructed as $0,1,2-$clique-sums of chordal graphs and/or cycles.
\QED

\section{$\text{CI}{\cal O}$ graphs and theta-ring graphs}\label{CIO-Theta}

In this section we prove that the $\text{CI}{\cal O}$ property is closed under $0,1,2-$clique-sums and that the chordal graphs are $\text{CI}{\cal O}$ graphs. In particular, we prove that if $D=G_{\cal O}$ is a tournament, then $P_D$ is a binomial complete intersection. Also we prove that prisms, pyramids, thetas and $\theta-$partial wheels (the forbidden induced subgraphs given in Corollary \ref{min-fobb-thetas}) are not $\text{CI}{\cal O}$ graphs. The main result in this section is the equivalence between $\text{CI}{\cal O}$ graphs and theta-ring graphs.

\begin{Definition}
A oriented complete graph is called a tournament.
\end{Definition}

\begin{Proposition}\label{completeCIO}\cite{Redei}
If $D_n=(K_n)_{{\cal O}_n}$ is a tournament, then $D_n$ has a Hamiltonian oriented path.
\end{Proposition}

\noindent
In the following result we will assume that: ${\cal O}$ is an orientation of $G$ and $x$ is a vertex of $G$ with $N_G(x)=\{x_1,\dotsc,x_{r_1}\}$ such that $P=(x_1,\dotsc,x_{r_1})$ is an oriented path in $D=G_{\cal O}$, i.e., $(x_j,x_{j+1}) \in E(D)$ for $j=1,\dotsc,r_1-1$.

\begin{Lemma}\label{pathorient2} 
Let ${\cal O}'$ be the orientation of $G'=G \setminus x$ induced by ${\cal O}$. If $C_i=(x,x_i,x_{i+1},x)$, then $P_D \subseteq (P_{D'} \cup \{t_{C_1},\dotsc,t_{C_{r_1-1}}\})$, where $D'={G'}_{{\cal O}'}$.
\end{Lemma}
\proof
For $1 \leq i < j \leq r_1$, we take the cycle $C_{i,j}=(x,x_i,P_{i,j},x_j,x)$ where $P_{i,j}$ is the oriented subpath of $P$ that joins $x_i$ and $x_j$. First we will prove that $t_{C_{i,j}} \in (t_{C_1},\dotsc,t_{C_{r_1-1}})$ by induction on $j-i$. If $j-i=1$ it is clear. Using the induction hypothesis we have that $t_{C_{i,j-1}} \in(t_{C_1},\dotsc,t_{C_{r_1-1}})$. Thus, by Lemma \ref{cyclegen} we have that $t_{C_{i,j}} \in (t_{C_{i,{j-1}}},t_{C_{j-1}})$. Therefore $t_{C_{i,j}} \in (t_{C_1},\dotsc,t_{C_{r_1-1}})$.
\\
Let $C$ be a cycle without chords of $D$. If $x \notin V(C)$ then $C \subset D'$ and $t_C \in P_{D'}$. Now, if 
$x \in V(C)$ then $C=(x,x_{j_1},L,x_{j_2},x)$ where $1 \leq j_1 < j_2 \leq r_1$ and $L$ is a path between $x_{j_1}$ and $x_{j_2}$ in $D$. Since, $C$ is a cycle without chords then, $V(L) \cap N_G(x)=\{x_{j_1},x_{j_2}\}$. On the other hand we take the cycle $C'=(x_{j_1}, L, x_{j_2}, P_{j_1,j_2},x_{j_1})$ then, $C' \subseteq D'$ and $t_{C'} \in P_{D'}$. By Lemma \ref{cyclegen} we have that $t_C \in (t_{C'},t_{C_{i,j}})$. Furthermore $t_{C_{i,j}} \in (t_{C_1}, \dotsc, t_{C_{{r_1}-1}})$ then, 
$t_C \in (P_{D'}, t_{C_1}, \dotsc, t_{C_{{r_1}-1}})$. Therefore 
$P_D \subseteq (P_{D'} \cup \{ t_{C_1},\dotsc,t_{C_{{r_1}-1}} \})$, by Proposition \ref{jul22-1-05}. 
\QED

\begin{Definition}
Let $G$ be a graph. A simplicial vertex $x$ is a vertex such that the induced subgraph obtained from $N_G(x)$ in $G$ is a complete graph. 
\end{Definition}

\begin{Proposition}\label{simplicial-chordal}\cite{dirac}
If $G$ is a chordal graph then, $G$ has a simplicial vertex. 
\end{Proposition}

\begin{Theorem}\label{Chordal-CIO}
If $G$ is a chordal graph then $G$ is a $\text{CI}{\cal O}$ graph. 
\end{Theorem}
\proof
By induction on $\vert V(G) \vert$. By Proposition \ref{componentsCIO} we can suppose that $G$ is a connected graph. Since $G$ is chordal there exists a vertex $x$ in $G$ such that $x$ is a simplicial vertex. Thus, if $K$ is the subgraph induced by $N_{G}(x)$, then $K$ is a complete graph and $G'=G \setminus x$ is a connected subgraph of $G$. We can suppose that $N_G(x)=\{x_1,\dotsc,x_s\}$. Let ${\cal O}$ be an orientation of $G$ and $D=G_{\cal O}$. By induction hypothesis $G'$ is $\text{CI}{\cal O}$, hence if ${\cal O}'$ is the orientation of $G'$ induced by ${\cal O}$ and $D'=G'_{{\cal O}'}$, then there exists the binomial generating set ${\cal G}'$ of $P_{D'}$  with $\vert {\cal G}' \vert=(q-s)-(n-1)+1$. On the other hand, by Proposition \ref{completeCIO}, there exists a Hamiltonian oriented path $P$ of $K$. We can suppose that $P=(x_1,\dotsc,x_s)$. We take the triangle $C_i=(x,x_i,x_{i+1},x)$ for $i=1,\dotsc, s-1$ and ${\cal G}={\cal G}' \cup \{t_{C_1},\dotsc,t_{C_{s-1}}\}$. Hence, by Lemma \ref{pathorient2}, $\cal G$ is a generating set of $P_D$ and $\vert {\cal G} \vert = q-n+1$. Therefore $G$ is $\text{CI}{\cal O}$.
\QED

\begin{Corollary}
If $D$ is a tournament, then $P_D$ is a binomial complete intersection.
\end{Corollary}
\proof
Since $D$ is a tournament then $D=K_{\cal O}$ where $K$ is a complete graph. Thus $K$ is a chordal graph and by Theorem~\ref{Chordal-CIO} $K$ is a $\text{CI}{\cal O}$ graph. Therefore $P_D$ is a binomial complete intersection.
\QED

\begin{Proposition}\label{2-sumCIO}
Let $G_1$, $G_2$ be two connected graphs and let $G$ be the $2-$clique-sum of $G_1$ and $G_2$. If $G_1$ and $G_2$ are $\text{CI}{\cal O}$ graphs then $G$ is a $\text{CI}{\cal O}$ graph.
\end{Proposition}
\proof
Let ${\cal O}$ be an orientation of the edges of $G$ and let ${\cal O}_i$ be the orientation of the edges of $G_i$ induced by ${\cal O}$, for $i=1,2$. If $\vert V(G_i) \vert =n_i$ and $\vert E(G_i) \vert =q_i$, for $i=1,2$, then there exists a set $B_i$ with $q_i-n_i+1$ binomials such that $B_i$ generates $P_{D_i}$, where $D_i=(G_i)_{{\cal O}_i}$. Since $G$ is the $2-$clique-sum of $G_1$ and $G_2$ then $\vert E(G) \vert = q_1+q_2-1$ and $\vert V(G) \vert = n_1+n_2-2$. If $C$ is a cycle without chords of $D$ then $V(C) \subset V(D_1)$ or $V(C) \subset V(D_2)$. Hence, $t_{C} \in (B_1 \cup B_2)$ and by Proposition \ref{jul22-1-05}, we have that $B_1 \cup B_2$ is a generating set of $P_{G_{\cal O}}$. But 
\[
\vert B_1 \cup B_2\vert=\sum_{i=1}^2 (q_i-n_i+1)=(q_1+q_2-1)-(n_1+n_2-2)+1= \mbox{ ht}(P_{G_{\cal O}}).
\]
Therefore $G$ is $\text{CI}{\cal O}$.
\QED

If $L=(t_{i_1},t_{i_2},\dotsc,t_{i_s})$ is an oriented path in $D=G_{\cal O}$, we denote by
$t_L$ the product of the edges of $L$, i.e., $t_L=t_{i_1}t_{i_2}\dotsm t_{i_s}$.

\begin{Proposition}\label{noCIO}
Prisms, pyramids, thetas and $\theta-$partial wheels are not $\text{CI}{\cal O}$ graphs.
\end{Proposition}
\proof
Case (a). $T$ is a partial wheel where $V(T)=\{x_1,\dotsc,x_m,x\}$ with center $x$, rim $C=(x_1,x_2,\dotsc,x_m,x_1)$, and $N_T(x)=\{x_{j_1},\dotsc,x_{j_k}\}$.
\\
First subcase, if $k \geq 4$ we take ${\cal T}=T_{\cal O}$ with acyclic orientation ${\cal O}$ given into Fi\-gure~\ref{O-pw4}.
\begin{figure}[h]
\setlength{\unitlength}{1mm}
\begin{picture}(50,49)
  \put(25,26){\circle*{2}}
  \put(18.6,45.1){\circle*{1.5}}
  \put(6,19.5){\circle*{1.5}}
  \put(31.5,7){\circle*{1.5}}
  \put(44,32.3){\circle*{1.5}}
  \put(31.3,45.1){\circle*{1.5}}
  \put(6,32.5){\circle*{1.5}}
  \put(18.5,7){\circle*{1.5}}
  \put(44,19.5){\circle*{1.5}}  
  \put(25,26){\vector(1,3){5}}
  \put(25,26){\line(1,3){6.35}}
  \put(25,26){\vector(1,-3){5}}
  \put(25,26){\line(1,-3){6.35}}
  \put(25,26){\vector(-3,1){15}}
  \put(25,26){\line(-3,1){19}}
  \put(25,26){\vector(-1,3){5}}
  \put(25,26){\line(-1,3){6.35}}
  \put(25,26){\vector(-3,-1){15}}
  \put(25,26){\line(-3,-1){19}}
  \put(25,26){\vector(-1,-3){5}}
  \put(25,26){\line(-1,-3){6.35}}
  \put(44,32.3){\vector(0,-1){7}}
  \put(44,32.3){\line(0,-1){12.8}}
  \put(6,19.5){\vector(0,1){9}}
  \put(6,19.5){\line(0,1){12.8}}
  \put(18.6,45.1){\vector(1,0){7}}
  \put(18.6,45.1){\line(1,0){12.8}}
  \put(18.6,45){\vector(-1,-1){8}}
  \put(18.6,45){\line(-1,-1){12.65}}
  \put(31.3,7){\vector(-1,0){7}}
  \put(31.3,7){\line(-1,0){12.6}}
  \put(6,19.5){\vector(1,-1){6}}
  \put(6,19.5){\line(1,-1){12.5}}
  \put(31.5,7){\line(1,1){12.5}}
  \put(44,32.3){\line(-1,1){12.8}}
  \put(44,19.5){\vector(-1,-1){7}}
  \put(31.3,45){\vector(1,-1){8}}
  \put(22,39){$t_{1}$}
  \put(30,37){$t_{2}$} 
  \put(30,14){$t_{k-3}$} 
  \put(15,16){$t_{k-2}$}
  \put(11,24){$t_{k-1}$}  
  \put(12,32){$t_{k}$} 
  \put(23,47){$L_{1}$} 
  \put(23,2){$L_{k-3}$} 
  \put(7,9){$L_{k-2}$} 	
  \put(0,25){$L_{k-1}$}  
  \put(6,39){$L_{k}$}
  \put(35,31){.}
  \put(35,26){.}
  \put(35,21){.}  
  \put(60,25){$
   \begin{array}{l}
   %\mbox{ The equations of the cycles without chords are: }\\
   g_1\quad = t_2-t_1t_{L_1}\\
   g_2\quad = t_3-t_2t_{L_2}\\
   \quad\quad\ \vdots \\
   g_i\quad = t_{i+1}-t_it_{L_i}\\
   \quad\quad\ \vdots \\
   g_{k-3} = t_{k-2}-t_{k-3}t_{L_{k-3}}\\
   g_{k-2} = t_{k-2}-t_{k-1}t_{L_{k-2}}\\
   g_{k-1} = t_k-t_{k-1}t_{L_{k-1}}\\
   g_k\quad = t_k-t_1t_{L_k}\\
   g_{k+1} = t_{L_{k-2}}t_{L_k}-t_{L_1}\dotsm t_{L_{k-3}}t_{L_{k-1}}
   \end{array}
  $}
\end{picture}
\caption{Acyclic orientation of the partial wheel ${\cal T}$ with $k\ge 4$, and the equations of its chordless cycles}\label{O-pw4}
\end{figure}

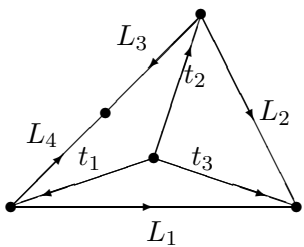
\begin{figure}[h]
\setlength{\unitlength}{1mm}
 \begin{picture}(60,32)
  \put(25,12){\circle*{1.5}}
  \put(6,5.5){\circle*{1.5}}
  \put(6,5.5){\line(1,1){25}}
  \put(6,5.5){\vector(1,1){7}}
  \put(31.3,31.1){\circle*{1.5}}
  \put(31.3,31.1){\vector(1,-2){7}}
  \put(31.3,31.1){\line(1,-2){13}}
  \put(31.3,31){\vector(-1,-1){7}}
  \put(18.6,18){\circle*{1.5}}
  \put(44,5.5){\circle*{1.5}}
  \put(6,5.5){\line(1,0){38}}
  \put(6,5.5){\vector(1,0){19}}  
  \put(25,12){\line(3,-1){19}}
  \put(25,12){\vector(3,-1){15}}
  \put(25,12){\vector(1,3){5}}
  \put(25,12){\line(1,3){6.35}}
  \put(25,12){\vector(-3,-1){15}}
  \put(25,12){\line(-3,-1){19}}
  \put(29,22){$t_2$} 
  \put(15,11){$t_1$}
  \put(30,11){$t_3$}
  \put(39,17){$L_2$}
  \put(24,1){$L_1$}
  \put(8,14){$L_4$}
  \put(20,27){$L_3$}       
  \put(60,17){$
     \begin{array}{l}
        %\mbox{The equations of the cycles without chords are:}\\
        g_1 = t_3-t_1t_{L_1}\\
        g_2 = t_3-t_2t_{L_2}\\
        g_3 = t_1t_{L_4}-t_2t_{L_3}\\
        g_4 = t_{L_1}t_{L_3}-t_{L_2}t_{L_4}
     \end{array}
  $}
 \end{picture}
 \caption{Acyclic orientation of the partial wheel ${\cal T}$ with $k=3$, and the equations of its chordless cycles}\label{O-pw3}
 \end{figure}

\noindent
We have that $\text{ht}(P_{\cal T})=k$. Let ${\cal G}=\{f_1,\dotsc,f_s\}$ be a minimum binomial set of generators for  $P_{\cal T}$. Since there are no oriented cycles in ${\cal T}$ and using Corollary \ref{oriencycle} then, there is no binomial having as a monomial 1. Since $g_i \in P_{\cal T}$ then, there exists a binomial $ f_{l_i} \in {\cal G}$ such that $f_{l_i}$ has as a monomial $t_i$ for $i=2,\dotsc,k-2$ and $i=k$. Without loss of generality we can suppose that $l_i=i$ for $i=2,\dotsc,k-2$ and $l_{k}=k-1$. Observe that if $C$ is a cycle of $T$ and $E(C) \cap E(L_i) \neq \emptyset$ then $E(L_i) \subset E(C)$. We will use this observation in the rest of the proof. 
\\
By the form of $g_{k+1}$ there exists a binomial $f_j \in {\cal G}$ with a monomial that divides $t_{L_{k-2}}t_{L_k}$. We can suppose that $f_j=t^{{\alpha}_j}-t^{{\beta}_j}$ and $t^{{\alpha}_j} \vert t_{L_{k-2}}t_{L_k}$. By Theorem \ref{bincycle} there exists a cycle $C_1$ such that the monomials of $t_{C_1}=t^{{\alpha}'_j}-t^{{\beta}'_j}$ divide the monomials of $f_j$. First we prove $E(L_k) \subset E(C_1)$. Suppose that it is not true, then $E(L_k) \cap E(C_1)=\emptyset$. Hence, $E(L_{k-2}) \cap E(C_1) \neq \emptyset$ and $E(L_{k-2}) \subset E(C_1)$. Furthermore since $t^{{\alpha}'_j} \vert t_{L_{k-2}}t_{L_k}$ then, $t_{k-1} \notin E(C_1)$ and $E(L_{k-1})\subset E(C_1)$. But $E(L_k) \cap E(C_1)=\emptyset$ then, $t_k \in E(C_1)$ and $t_k \vert t^{{\alpha}'_j}$. This is a contradiction because $t_k$ does not divide $t_{L_{k-2}}t_{L_k}$. Therefore $E(L_k) \subset E(C_1)$. Now, we prove that $x \notin E(C_1)$. Suppose $x \in E(C_1)$, we take the minimal $i_1=\min\{ i \in \{1,\dotsc,k\} \mid t_i \in E(C_1) \}$. Thus, $\cup^{i_1}_{i=1} E(L_i) \subset E(C_1)$ and $t_{i_1} \vert t^{{\alpha}'_j}$. But this is not possible because $t_{i_1}$ does not divide $t_{L_{k-2}}t_{L_k}$. Hence, $x \notin E(C_1)$ and $C_1={\cal T} \setminus x$. Therefore $t_{C_1}=g_{k+1}$. Using Proposition \ref{nwell-div} we have that $f_j=g_{k+1}$. Furthermore $f_j \neq f_i$ for $i=2,\dotsc,k-1$, then we can suppose that $j=1$.
\\
On the other hand by the form of $g_{k-1}$ there exists a binomial $f_{j_1} \in {\cal G}$ with a monomial 
that divides $t_{k-1}t_{L_{k-1}}$. By Theorem \ref{bincycle} there exists a cycle $C_2$ such that the monomials
of $t_{C_2}=t^{{\alpha}'_{j_1}}-t^{{\beta}'_{j_1}}$ divide the monomials of $f_j$. We can suppose that 
$t^{{\alpha}'_{j_1}}\vert t_{k-1}t_{L_{k-1}}$. Suppose $E(L_{k-1})\cap E(C_2)=\emptyset$ then $t_{k-1} \in E(C_2)$ and $E(L_{k-2}) \subset E(C_2)$. Thus, $t_{L_{k-2}} \vert t^{{\alpha}'_{j_1}}$ which is not possible because $t_{L_{k-2}}$ does not divide $t_{k-1}t_{L_{k-1}}$. So, $E(L_{k-1})\cap E(C_2) \neq \emptyset$ and $E(L_{k-1})\subset E(C_2)$. Now, we prove that $t_{k-1} \in E(C_2)$. If $t_{k-1} \notin E(C_2)$ then $E(L_{k-2}) \subset E(C_2)$. Furthermore, $E(L_{k-3}) \subset E(C_2)$ or $t_{k-2} \in E(C_2)$. Hence, $t_{L_{k-3}} \vert t^{{\alpha}'_{j_1}}$ or $t_{k-2} \vert t^{{\alpha}'_{j_1}}$. But this is not possible, then $t_{k-1} \in E(C_2)$. Now, $t_k \in E(C_2)$ or $E(L_k) \subset E(C_2)$. If $t_k \in E(C_2)$ then $C_2=(t_{k-1},L_{k-1},t_k)$ and $t_{C_2}=-g_{k-1}$. For the other case, if $E(L_k) \subset E(C_2)$ then $C_2=(t_{k-1},L_{k-1},L_{k},t_1)$ and $t_{C_2}=t_{k-1}t_{L_{k-1}}-t_1t_{L_{k}}=g_{k}-g_{k-1}$. In both cases, by Proposition \ref{nwell-div} we have that $f_{j_1}=t_{C_2}$. Furthermore, by the form of $g_k$ there exists a binomial $f_{j_2} \in {\cal G}$ with a monomial that divides $t_1t_{L_k}$. In a similar way as in the last argument, we obtain that there exists a cycle $C_3$ such that $f_{j_2}=t_{C_3}$ where $C_3=(t_k,L_k,t_1)$ or $C_3=(t_{k-1},L_{k-1},L_{k},t_1)$. So, $f_{j_1}, f_{j_2} \in \{t_k-t_{k-1}t_{L_{k-1}}, t_k-t_1t_{L_k}, t_{k-1}t_{L_{k-1}}-t_1t_{L_{k}} \}$. If $\{f_1,\dots,f_{k-1}\} \cap \{f_{j_1},f_{j_2}\} \neq \emptyset$ then $f_{j_1}=f_{k-1}$ or $f_{j_2}=f_{k-1}$ and in this case $f_{j_1} \neq f_{j_2}$. Therefore $k \leq \vert\{f_1,\dotsc,f_{k-1},f_{j_1},f_{j_2}\}\vert$.
\\
By the form of $g_{k-2}$ there exists a binomial $f_{j_3}$ with a monomial that divides $t_{k-1}t_{L_{k-2}}$.   
Thus, there exists a cycle $C_4$ such that the monomials of $t_{C_4}=t^{{\alpha}'_{j_3}}-t^{{\beta}'_{j_3}}$ divide the monomials of $f_{j_3}$. We can suppose that $t^{{\alpha}'_{j_3}}\vert t_{k-1}t_{L_{k-2}}$. Since $t^{{\alpha}'_{j_3}} \neq 1$ and $t_{L_{k-1}}$ does not divide $t^{{\alpha}'_{j_3}}$ then $V(L_{k-2}) \subset V(C_4)$. Furthermore, since $t_{L_k}$ and $t_k$ do not divide $t^{{\alpha}'_{j_3}}$ then $t_{k-1} \in E(C_4)$.
Hence, if $t_{k-2} \in E(C_4)$ then $C_4=(t_{k-2},L_{k-2},t_{k-1})$ and $t_{C_4}=t_{k-1}t_{L_{k-2}}-t_{k-2}$.
In the other case, if $t_{k-2} \notin E(C_4)$ then $t_{C_4}=t_{k-1}t_{L_{k-2}}-t^{{\beta}'_{j_3}}$. Thus, $t^{{\alpha}'_{j_3}}=t_{k-1}t_{L_{k-2}}$ and by Proposition \ref{nwell-div} we have that $f_{j_3}=t_{C_4}$. Furthermore, if $t_{k-2} \notin E(C_4)$ then $k+1 \leq \vert\{f_1,\dotsc,f_{k-1},f_{j_1},f_{j_2},f_{j_3}\}\vert$ and $P_{\cal T}$ is not a binomial complete intersection. So, we can assume that $t_{k-2} \in E(C_4)$ and $f_{j_3}=t_{k-1}t_{L_{k-2}} - t_{k-2}$. Moreover, if $f_{j_3} \neq -f_{k-2}$ then $k+1 \leq \vert\{f_1,\dotsc,f_{k-1},f_{j_1},f_{j_2},f_{j_3}\}\vert$ and $P_{\cal T}$ is not a binomial complete intersection. Hence, we can assume $-f_{k-2}=f_{j_3}=t_{k-1}t_{L_{k-2}} - t_{k-2}$ and it is the unique binomial in ${\cal G}$ which has a monomial that divides either $t_{k-2}$ or $t_{k-1}t_{L_{k-2}}$. Then,
\[
\begin{array}{ccc}
t_{k-2} - t_{k-3} t_{L_{k-3}}=\sum_{i=1}^s{g_if_i} & = & g_{k-2}(t_{k-1}t_{L_{k-2}}-t_{k-2})+\sum_{i\neq k-2}{g_i f_i}\\
& = & (-1+g'_{k-2})(t_{k-1}t_{L_{k-2}}-t_{k-2})+\sum_{i \neq k-2}{g_i f_i}.
\end{array}
\]
Where $\deg(g'_{k-2})>0$ or $g'_{k-2}=0$. Using the last equation we obtain  
\[
t_{k-1}t_{L_{k-2}}- t_{k-3}t_{L_{k-3}}=g'_{k-2}(t_{k-1}t_{L_{k-2}}-t_{k-2})+\sum_{i \neq k-2}{g_if_i}.
\]
So, there exists a $f_i$ with $i \neq k-2$ where one of its monomials divides $t_{k-1}t_{L_{k-2}}$. But it is not possible.
\\
\\
Second subcase, if $k=3$ we take ${\cal T}=T_{\cal O}$ with the acyclic orientation ${\cal O}$ given into Figure~\ref{O-pw3}.

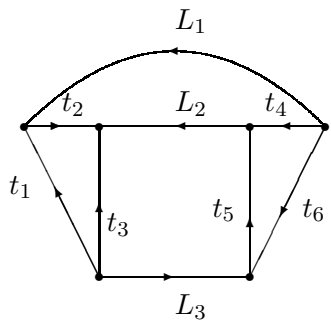
\begin{figure}[h]
\setlength{\unitlength}{1mm}
\begin{picture}(50,42)
  \put(45,27){\circle*{1}}
  \put(35,27){\circle*{1}}
  \put(35,7){\circle*{1}}
  \put(15,27){\circle*{1}}
  \put(15,7){\circle*{1}}
  \put(5,27){\circle*{1}}
  \put(45,27){\vector(-1,0){6}}
  \put(45,27){\line(-1,0){9.5}}
  \put(35,7){\vector(0,1){8}}
  \put(35,7){\line(0,1){20}}
  \put(45,27){\vector(-1,-2){6}}
  \put(45,27){\line(-1,-2){10}}
  \put(35,27){\vector(-1,0){10}}
  \put(35,27){\line(-1,0){20}}
  \put(15,7){\vector(1,0){10}}
  \put(15,7){\line(1,0){20}}
  \put(15,7){\vector(0,1){10}}
  \put(15,7){\line(0,1){20}}
  \put(15,7){\vector(-1,2){6}}
  \put(15,7){\line(-1,2){10}}
  \put(5,27){\vector(1,0){5}}
  \put(5,27){\line(1,0){10}}
  \put(5,27){\vector(1,0){5}}
  \put(5,27){\line(1,0){10}}
  \qbezier(5,27)(25,47)(45,27)
  \put(24.5,37){\vector(-1,0){.5}}
  \put(10,29){$t_{2}$}
  \put(37,29){$t_{4}$} 
  \put(3,18){$t_{1}$} 
  \put(16,13){$t_{3}$}  
  \put(30,15){$t_{5}$} 
  \put(42,15){$t_{6}$}
  \put(25,40){$L_1$}
  \put(25,29){$L_2$}
  \put(25,2){$L_3$}
  \put(60,25){$
   \begin{array}{l}
   g_1 = t_4-t_5t_6\\
   g_2 = t_3-t_1t_2\\
   g_3 = t_3-t_{L_3}t_5t_{L_2}\\
   g_4 = t_1t_6-t_{L_1}t_{L_3}\\
   g_5 = t_4t_{L_2}-t_2t_{L_1}
   \end{array}
   $}
 \end{picture}
 \caption{Acyclic orientation of the prism ${\cal T}$, and the equations of its cycles without chords}\label{O-prism}
 \end{figure}

\noindent
We have that $\text{ht}(P_{\cal T})=3$. Let ${\cal G}=\{f_1,\dotsc,f_s\}$ be a minimum binomial set of 
generators of  $P_{\cal T}$. For $m \in A=\{t_3, t_1t_{L_1},t_2t_{L_2}, t_1t_{L_4}, t_2t_{L_3}, t_{L_1}t_{L_3}, t_{L_2}t_{L_4}\}$ there exists a binomial $f \in {\cal G}$ such that $f$ has a monomial $t^{\alpha}$ and 
$t^{\alpha} \vert m$. By Theorem \ref{bincycle} there exists a cycle $C$ such that $t_{C}=t^{{\alpha}'}-t^{{\beta}'}$ and $t^{{\alpha}'} \vert t^{\alpha}$. But there are no cycles whose binomials have either $t_{L_i}$ or $t_j$ as a monomial, for $i=1,\dotsc,4$ and $j=1,2$. Hence, $t^{{\alpha}'}=t^{\alpha}=m$ and by Proposition \ref{nwell-div} $t_C=f$. Then, for every $m \in A$ there exists $f \in {\cal G}$ such that $m$ is a monomial of $f$. But $\vert A \vert/2 \leq \vert {\cal G} \vert$ and $\vert A \vert=7$ then $4 \leq \vert {\cal G} \vert$. Therefore $P_{\cal T}$ is not a binomial complete intersection.   
\\
\\
Case (b-1). If $T$ is a prism, we take ${\cal T}=T_{\cal O}$ with the acyclic orientation ${\cal O}$ given into Fi\-gure~\ref{O-prism}.

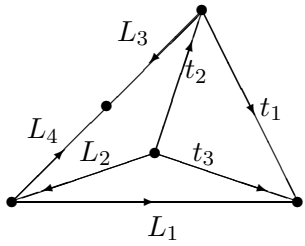
\begin{figure}[h]
 \setlength{\unitlength}{1mm}
 \begin{picture}(50,32)
  \put(25,12){\circle*{1.5}}
  \put(6,5.5){\circle*{1.5}}
  \put(6,5.5){\line(1,1){25}}
  \put(6,5.5){\vector(1,1){7}}
  \put(31.3,31.1){\circle*{1.5}}
  \put(31.3,31.1){\vector(1,-2){7}}
  \put(31.3,31.1){\line(1,-2){13}}
  \put(31.3,31){\vector(-1,-1){7}}
  \put(18.6,18.3){\circle*{1.5}}
  \put(44,5.5){\circle*{1.5}}
  \put(6,5.5){\line(1,0){38}}
  \put(6,5.5){\vector(1,0){19}}  
  \put(25,12){\line(3,-1){19}}
  \put(25,12){\vector(3,-1){15}}
  \put(25,12){\vector(1,3){5}}
  \put(25,12){\line(1,3){6.35}}
  \put(25,12){\vector(-3,-1){15}}
  \put(25,12){\line(-3,-1){19}}
  \put(29,22){$t_2$} 
  \put(15,11){$L_2$}
  \put(30,11){$t_3$}
  \put(39,17){$t_{1}$}
  \put(24,1){$L_1$}
  \put(8,14){$L_4$}
  \put(20,27){$L_3$}       
  \put(60,17){$\begin{array}{l}
        %\mbox{The equations of the cycles without chords are:}\\
        g_1 = t_3-t_{L_1}t_{L_2}\\
        g_2 = t_3-t_1t_{2}\\
        g_3 = t_{L_2}t_{L_4}-t_2t_{L_3}\\
        g_4 = t_{L_1}t_{L_3}-t_{1}t_{L_4}
     \end{array}$}
 \end{picture}
 \caption{Acyclic orientation of the pyramid ${\cal T}$, and the equations of its cycles without chords}\label{O-pyramid}
\end{figure}

\noindent
We have that $\text{ht}(P_{\cal T})=4$. Let ${\cal G}=\{f_1,\dotsc,f_s\}$ be a minimum binomial set of generators of  $P_{\cal T}$. We take $A=\{t_5t_6, t_1t_2, t_2t_{L_1}, t_{L_3}t_5t_{L_2}, t_1t_6, t_{L_1}t_{L_3}, t_3, t_4\}$. By the same argument as for partial wheels with $k=3$, we have that for every
$m \in A$ there exists $f \in {\cal G}$ and a cycle $C$ of $D$ such that $m$ is a monomial of $f$
and $f=t_{C}$. Let $f_{j_1}$, $f_{j_2}$, $f_{j_3}$ be binomials of ${\cal G}$ that have as a monomial $t_3$, $t_1t_2$ and $t_{L_3}t_5t_{L_2}$, respectively. Hence, if $f_{j_1}=t_{C_1}$, $f_{j_2}=t_{C_2}$ and 
$f_{j_3}=t_{C_3}$ then, 
\[
C_1,C_2,C_3 \in \{(t_3, t_1, t_2), (t_1, t_2, L_2, t_5, L_3), (t_3, L_3, t_5, L_2)\}.
\] 
Thus, if $B=\{f_{j_1},f_{j_2},f_{j_3}\}$ then $\vert B \vert \geq 2$. Since if $f_{j_1}=f_{j_3}$ then $C_2=C_3=(t_1, t_2, L_2, t_5, L_3)$ and $f_{j_1} \neq f_{j_2}$. Furthermore for every monomial in $A'=A \setminus \{t_3, t_1t_2, t_{L_3}t_5t_{L_2}\}$ there exists a binomial in ${\cal G} \setminus B$. But 
$\vert {\cal G} \setminus B \vert \geq \vert A' \vert/2$ and $\vert A' \vert=5$ then $\vert {\cal G} \setminus B \vert \geq 3$ and $\vert {\cal G} \vert \geq 5$. Therefore $P_{\cal T}$ is not a binomial complete intersection.   
\\
\\
Case (b-2). If $T$ is a pyramid, we take ${\cal T}=T_{\cal O}$ with the acyclic orientation ${\cal O}$ given into Figure~\ref{O-pyramid}.

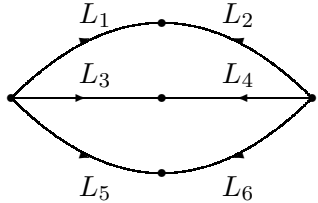
\begin{figure}[h]
 \setlength{\unitlength}{1mm}
 \begin{picture}(50,27)
  \put(45,15){\circle*{1}}
  \put(25,15){\circle*{1}}
  \put(5,15){\circle*{1}}
  \put(45,15){\vector(-1,0){10}}
  \put(45,15){\line(-1,0){40}}
  \put(5,15){\vector(1,0){10}}
  \put(5,15){\line(1,0){10}}
  \qbezier(5,15)(25,35)(45,15)
  \qbezier(5,15)(25,-5)(45,15)
  \put(25,25){\circle*{1}}
  \put(25,5){\circle*{1}}
  \put(14,17){$L_3$}
  \put(33,17){$L_4$}
  \put(15,22.8){\vector(3,1){.5}}
  \put(35,22.8){\vector(-3,1){.5}} 
  \put(14,25){$L_1$}
  \put(33,25){$L_2$}
  \put(15,7.2){\vector(3,-1){.5}}
  \put(35,7.2){\vector(-3,-1){.5}} 
  \put(14,2){$L_5$}
  \put(33,2){$L_6$}
  \put(60,14){$
     \begin{array}{l}
        %\mbox{The equations of the cycles without chords are:} \\
        g_1 = t_{L_1}t_{L_4}-t_{L_2}t_{L_3}\\
        g_2 = t_{L_3}t_{L_6}-t_{L_4}t_{L_5}\\
        g_3 = t_{L_1}t_{L_6}-t_{L_2}t_{L_5}
     \end{array}$}
 \end{picture}
 \caption{Acyclic orientation of the theta ${\cal T}$, and the equations of its cycles without chords}\label{O-theta}
\end{figure}

\noindent 
This case is similar to the case when $T$ is a partial wheel with $k=3$, only change $t_1$ by the path $L_2$. Therefore $P_{\cal T}$ is not a binomial complete intersection. 
\\
\\
Case (b-3) If $T$ is a theta, we take ${\cal T}=T_{\cal O}$ with the acyclic orientation ${\cal O}$ given into Figure~\ref{O-theta}. 

\noindent
We have that $\text{ht}(P_{\cal T})=2$. Let ${\cal G}=\{f_1,\dotsc,f_s\}$ be a minimum binomial set of 
generators of  $P_{\cal T}$. We take $A=\{t_{L_1}t_{L_4}, t_{L_2}t_{L_3}, t_{L_3}t_{L_6}, t_{L_4}t_{L_5}, t_{L_1}t_{L_6}, t_{L_2}t_{L_5}\}$. By the same argument as for partial wheels with $k=3$, we have that for every $m \in A$ there exists $f \in {\cal G}$ and a cycle $C$ of $D$ such that $m$ is a monomial of $f$ and $f=t_{C}$. Hence, $\vert {\cal G} \vert \geq \vert A \vert/2$ but $\vert A \vert=6$ then $\vert {\cal G} \vert \geq 3$. Therefore $P_{\cal T}$ is not a binomial complete intersection.
\QED

\begin{Theorem}\label{For-all-Theta-exists-deltaCIO}
Let $G$ be a graph. $G$ is a $\text{CI}{\cal O}$ graph if and only if $G$ is a theta-ring graph.
\end{Theorem}
\proof
$\Rightarrow$) If $G$ is not a theta-ring graph then, by Theorem~\ref{For-all-Theta-exists-delta}, $G$ has an induced subgraph $H$ that is a prism, pyramid, theta or $\theta-$partial wheel. Thus, by Proposition~\ref{noCIO}, $H$ is not a $\text{CI}{\cal O}$ graph. Hence, by Corollary \ref{inducedCIO}, $G$ is not a $\text{CI}{\cal O}$ graph which is a contradiction. Therefore $G$ is a theta-ring graph.
\\
$\Leftarrow$) By Theorem \ref{For-all-Theta-exists-delta} $G$ can be constructed by $0,1,2-$clique-sums of chordal graphs and/or cycles. By Proposition \ref{Chordal-CIO} chordal graphs and cycles are $\text{CI}{\cal O}$ graphs. Therefore by Proposition \ref{componentsCIO}, Proposition \ref{block-CIO} and Proposition \ref{2-sumCIO}, $G$ is a $\text{CI}{\cal O}$ graph.
\QED

\begin{Definition}\cite{GRV}
A graph $G$ is a {\it ring graph\/} if each block of $G$ which is not a bridge or a vertex can be constructed successively by $2-$clique-sums of cycles. 
\end{Definition}

\begin{Corollary}\cite{Ring}
If $G$ is a ring graph, then $G$ is a $\text{CI}{\cal O}$ graph. The converse holds if $G$ is bipartite. 
\end{Corollary} 
\proof
First we observe that ring graphs can be obtained by $0,1,2-$clique-sums of vertices, edges and cycles. 
Hence, by Theorem \ref{For-all-Theta-exists-delta}, $G$ is a theta-ring graph and by Theorem \ref{For-all-Theta-exists-deltaCIO}, $G$ is a $\text{CI}{\cal O}$ graph.
\\
Now, if $G$ is bipartite and $G$ is a $\text{CI}{\cal O}$ graph, then by Theorem \ref{For-all-Theta-exists-delta} and Theorem \ref{For-all-Theta-exists-deltaCIO}, $G$ can be constructed by $0,1,2-$clique-sums of chordal graphs and/or cycles. Since $G$ is bipartite, if $H$ is a chordal induced subgraph of $G$, then $H$ is a forest. Hence, $G$ can be constructed by $0,1,2-$clique-sums of cycles and/or edges. Therefore, $G$ is a ring graph.
\QED

\bigskip

{\bf Acknowledgements}  The authors thank ABACUS, CONACyT grant EDOMEX-2011-C01-165873.

\bibliographystyle{plain}

\end{document}